\documentclass[a4paper,12pt]{amsart}
\usepackage{amssymb}
\input xypic
\xyoption{all}
\def\hbar{\bar{h}}

\def\iso{\buildrel \sim\over\to}

\def\GS{{\mathfrak{S}}}

\def\CC{{\mathcal{C}}}

\def\CF{{\mathcal{F}}}
\def\CG{{\mathcal{G}}}

\def\CL{{\mathcal{L}}}

\def\CO{{\mathcal{O}}}

\def\CR{{\mathcal{R}}}

\def\BA{{\mathbf{A}}}
\def\BB{{\mathbf{B}}}

\def\BF{{\mathbf{F}}}
\def\BG{{\mathbf{G}}}

\def\BL{{\mathbf{L}}}

\def\BP{{\mathbf{P}}}

\def\BT{{\mathbf{T}}}
\def\BU{{\mathbf{U}}}
\def\BV{{\mathbf{V}}}

\def\BZ{{\mathbf{Z}}}

\def\Br{{\mathbf{r}}}

\def\br{\operatorname{br}\nolimits}
\def\Br{\operatorname{Br}\nolimits}

\def\can{{\mathrm{can}}}

\def\End{\operatorname{End}\nolimits}

\def\GL{\operatorname{GL}\nolimits}

\def\Hom{\operatorname{Hom}\nolimits}

\def\Ind{\operatorname{Ind}\nolimits}

\def\Irr{\operatorname{Irr}\nolimits}

\def\mMod{\operatorname{\!-mod}\nolimits}

\def\mMOD{\operatorname{\!-Mod}\nolimits}

\def\Res{\operatorname{Res}\nolimits}

\def\SL{\operatorname{SL}\nolimits}

\def\ie{{\em i.e.}}
\def\eg{{\em e.g.}}
\def\Qlbar{{\bar{\mathbf{Q}}_l}}

\def\oppose{{\mathrm{opp}}}

\newtheorem{thm}{Theorem}[section]
\newtheorem{lemma}[thm]{Lemma}
\newtheorem{cor}[thm]{Corollary}
\newtheorem{prop}[thm]{Proposition}

\newtheorem{conj}[thm]{Conjecture}

\theoremstyle{definition}
\newtheorem{rem}[thm]{Remark}

\begin{document}

\title{Coxeter orbits and modular representations}
\author{C\'edric Bonnaf\'e and Rapha\"el Rouquier}
\date{\today}
\address{C\'edric Bonnaf\'e~:
        Univ. de Franche-Comt\'e, D\'epartement de Math\'ematiques 
	(CNRS UMR 6623), 16 Route de Gray, 25000 Besan\c{c}on, FRANCE.}
\email{bonnafe@math.univ-fcomte.fr}
\address{Rapha\"el Rouquier\\
Institut de Math\'ematiques de Jussieu --- CNRS\\
UFR de Math\'ematiques, Universit\'e Denis Diderot\\
2, place Jussieu\\
75005 Paris, FRANCE\\
and
Department of Pure Mathematics, University of Leeds, Leeds LS2 9JT, U.K.}
\email{rouquier@math.jussieu.fr}

\begin{abstract}
We study the modular representations of finite groups of Lie type arising
in the cohomology of certain quotients of Deligne-Lusztig varieties associated
with Coxeter elements. These quotients are related to Gelfand-Graev
representations and we present a conjecture on the Deligne-Lusztig restriction
of Gelfand-Graev representations. We prove the conjecture for restriction to
a Coxeter torus. We deduce a proof of Brou\'e's conjecture on equivalences
of derived categories arising from Deligne-Lusztig varieties, for
a split group of type $A_n$ and a Coxeter element. Our study is based on
Lusztig's work in characteristic $0$ \cite{Lu2}.
\end{abstract}

\maketitle
\section*{Introduction}

In \cite{Lu2}, Lusztig proved that the Frobenius eigenspaces
on the $\Qlbar$-cohomology spaces of the Deligne-Lusztig variety $X$
associated to a Coxeter element are irreducible unipotent representations.
We give here a partial modular analog of this result. Except in type $A$,
we treat only part of the representations (those occurring in Harish-Chandra
induced Gelfand-Graev representations). Note that we need to consider 
non-unipotent representations
at the same time. We show Brou\'e's conjecture on
derived equivalences coming from Deligne-Lusztig varieties, in type $A$, and
for Coxeter elements (Theorem \ref{ADC}). These are the first non-trivial
examples with varieties of dimension $\ge 2$.

\medskip
Recall that Brou\'e's abelian defect group conjecture \cite{Br1}
predicts that an $\ell$-block of a finite group $G$ is derived
equivalent to the corresponding
block of the normalizer $H$ of a defect group $D$, when $D$ is abelian.
When the group under consideration is a finite group of Lie type and $\ell$
is not the defining characteristic (and not too small), 
then Brou\'e further conjectures
\cite{Br1} that the complex of cohomology of a certain Deligne-Lusztig
variety $Y$ should provide a complex realising an equivalence.
The variety $Y$ has a natural action of $G\times C_G(D)^\oppose$, which does
not extend in general to an action of $G\times H^\oppose$. Nevertheless, instead
of an action of the relative Weyl group $H/C_G(D)$, it is expected \cite{Br2}
that there
will be an action of its associated braid monoid, which will induce an action
of the corresponding Hecke algebra in $\End_{D^b(\BZ_{\ell}G)}(C)$, where
$C=R\Gamma_c(Y,\BZ_\ell)$ is the complex of proper support cohomology of $Y$.
Finally, the Hecke algebra should be isomorphic to the group algebra of
the relative Weyl group.
When $H/C_G(D)$ is cyclic, the required action should be provided by the
Frobenius endomorphism $F$.

In \cite{Rou2}, the case where $Y$ is a curve
was dealt with, and the key point was a (modular version of the)
disjunction property of cohomology spaces for the $G$-action, which was
a consequence of specific properties of curves.

\smallskip
Here, the key step is a disjunction property for the action of $T\rtimes F$.
Our approach uses very little information on the modular representations of
$G$ and might apply to other situations.
The crucial geometrical part is an explicit description of the quotient
of the Deligne-Lusztig variety $Y$ by $D(\BU_0)^F$ (\S \ref{secquo}),
the rational points
of the derived subgroup of the unipotent radical of an $F$-stable Borel
subgroup (Lusztig had shown that $\BU_0^F\backslash Y/T$ is a product of
$\BG_m$'s). We deduce a disjunction property
for the action of $T\rtimes F$ on the complex of cohomology of
$D(\BU_0)^F\backslash Y$ (Corollary \ref{specE}).

Considering the complex of cohomology of $D(\BU_0)^F\backslash Y$ amounts
to tensoring the complex of $Y$ by the sum of the Harish-Chandra
induced Gelfand-Graev representations. For $\GL_n$, such a sum is a
progenerator and this gives the expected disjunction property for $Y$.

We conjecture (\S \ref{secconjectures})
that the Deligne-Lusztig restriction of a Gelfand-Graev
representation is a shifted Gelfand-Graev representation. Such a result
is known at the level of $K_0$ and in the case of Harish-Chandra induction.
We deduce the truth of the conjecture for Coxeter tori, as a consequence
of our geometric study (Theorem \ref{coxeter}).

\medskip
This is our second paper (after \cite{BonRou}) attempting to
extend some of the fundamental results of Lusztig to the modular setting.
We dedicate this paper to Professor Lusztig on his sixtieth birthday.

\tableofcontents

\bigskip

\section{Preliminaries}
The results of this chapter are mostly classical.
\subsection{Definitions}
\subsubsection{}
Let $R$ be a commutative ring. Given $d$ a positive integer,
we denote by $\mu_d(R)$ the subgroup of
$R^\times$ of elements of order dividing $d$.

Given a set $I$ and a group $G$ acting on $I$, we denote by $[G\backslash I]$ a
subset of $I$ of representatives of orbits.

Let $A$ be an $R$-algebra. We denote by $A^\oppose$ the opposite algebra to $A$.
Similarly, if $G$ is a group, we denote by $G^\oppose$ the opposite group to
$G$ and we put $G^\#=G-\{1\}$.
Given $M$ an $R$-module, we put $AM=A\otimes_R M$.
We denote by $A\mMOD$ the category of $A$-modules and by
$A\mMod$ the category of finitely generated $A$-modules.
Given $\CC$ an additive category, we denote by
$K^b(\CC)$ the homotopy category of bounded complexes of objects of $\CC$.
When $\CC$ is an abelian category, 
we denote by $D^b(\CC)$ its bounded derived category.
We put $K^b(A)=K^b(A\mMod)$ and $D^b(A)=D^b(A\mMod)$ (when $A\mMod$ is an
abelian category).

Given $C$ and $D$ two complexes of $A$-modules, we denote by
$\Hom^\bullet_A(C,D)$ the total complex of the double complex
$(\Hom_A(C^i,D^j))_{i,j}$ and by
$R\Hom^\bullet_A(C,D)$ the corresponding derived version.

\subsubsection{}
Let $\ell$ be a prime number. Let $K$ be a finite extension
of the field of $\ell$-adic numbers, let $\CO$ be the normal
closure of the ring of
$\ell$-adic integers in $K$ and let $k$ be the residue field of
$\CO$. We assume $K$ is big enough for the finite groups under consideration
(\ie, $K$ contains the $e$-th roots of unity where $e$ is the exponent of
one of the finite groups considered).

Let $H$ be a finite group. We denote by $\Irr(H)$ the set of irreducible
characters of $H$ with values in $K$ and we put
$\Irr(H)^\#=\Irr(H)-\{1\}$. Given $\chi\in\Irr(H)$, we put
$e_\chi=\frac{\chi(1)}{|H|}\sum_{h\in H}\chi(h)h^{-1}$. If $\chi(1)=1$
and $R=K$, $\CO$ or $k$,
we denote by $R_\chi$ the $RH$-module $R$ on which
$H$ acts via $\chi$.

\subsubsection{}
Let $p$ be a prime number distinct from $\ell$.
We denote by $\bar{\BF}_p$
an algebraic closure of $\BF_p$. Given $q$ a power of $p$, we denote
by $\BF_q$ the subfield of $\bar{\BF}_p$ with $q$ elements.
Given $d$ a positive integer, we put $\mu_d=\mu_d(\bar{\BF}_p)$.

\subsubsection{}
Let $X$ be an algebraic variety over $\bar{\BF}_p$ acted on by
a finite group $H$. Let $R$ be a ring amongst $K$, $\CO$ and $k$. We denote
by $\tilde{R}\Gamma_c(X,R)$ the object of $K^b(R H\mMOD)$ representing
the cohomology with compact support
of $X$ in $R$, as defined in \cite{Ri1,Rou2}. It is
a bounded complex of $R H$-modules which are direct summands of
permutation modules (not finitely generated). If $X$ is equipped with
a Frobenius endomorphism, then that endomorphism induces an invertible
operator of $\tilde{R}\Gamma_c(X,R)$. Note finally that, as a complex of
$RH$-modules, $\tilde{R}\Gamma_c(X,R)$ is homotopy equivalent to
a bounded complex of finitely generated modules which are direct summands
of permutation modules.

We will denote by $R\Gamma_c(X,R)$ the corresponding ``classical''
object of $D^b(RH\mMOD)$. If $X$ is equipped with
a Frobenius endomorphism, then, $R\Gamma_c(X,R)$, viewed as a complex
with an action of $H$ and of the Frobenius endomorphism, is quasi-isomorphic
to a bounded complex of modules which have finite rank over $R$.

\subsection{Algebraic groups}
\subsubsection{}
Let $\BG$ be a connected reductive algebraic group over $\bar{\BF}_p$, with
an endomorphism $F$. We assume there is a positive integer $f$ such that
$F^f$ is a Frobenius endomorphism of $\BG$.

Let $\BB_0$ be an $F$-stable
Borel subgroup of $\BG$, let $\BT_0$ be an $F$-stable maximal torus of $\BB_0$,
let $\BU_0$ denote the unipotent radical of $\BB_0$, and let
$W=N_\BG(\BT_0)/\BT_0$
denote the Weyl group of $\BG$ relative to $\BT_0$.
Given $w \in W$, we denote by $\dot{w}$ a
representative of $w$ in $N_{\BG}(\BT_0)$. If moreover $w \in W^{F^m}$ for some
positive integer $m$,
then $\dot{w}$ is chosen in $N_{\BG^{F^m}}(\BT_0)$.

Let $\Phi$ denote the root system
of $\BG$ relative to $\BT_0$, $\Phi^+$ the set of positive roots of $\Phi$
corresponding to $\BB_0$ and $\Delta$ the basis of $\Phi$ contained in $\Phi^+$.
We denote by $\phi : \Phi \to \Phi$ the bijection such that $F(\alpha)$ is
a positive multiple of $\phi(\alpha)$. We denote by $d$ the order
of $\phi$. Note that $F^d$ is a Frobenius endomorphism of $\BG$ defining
a split
structure over a finite field with $q^d$ elements, where $q$ is a positive
real number. 

If $\alpha \in \Phi$, we denote by $s_\alpha$
the reflection with respect to $\alpha$, by $\alpha^\vee$ the associated
coroot, by
$q_\alpha^\circ$ the power of $p$ such that
$F(\alpha)=q_\alpha^\circ \phi(\alpha)$ and
by $d_\alpha$ the smallest natural
number such that $F^{d_\alpha}(\alpha)$ is a multiple of $\alpha$.
In other words,
$d_\alpha$ is the length of the orbit of $\alpha$ under the action of $\phi$.
We have $d={\mathrm{lcm}}(\{d_\alpha\}_{\alpha\in\Phi})$.
Let $q_\alpha=q_\alpha^\circ q_{\phi(\alpha)}^\circ \cdots
 q_{\phi^{d_\alpha-1}(\alpha)}^\circ$. We have
$q_\alpha=q^{d_\alpha}$.
Note that $d_{\phi(\alpha)}=d_\alpha$ and $q_{\phi(\alpha)}=q_\alpha$.
Let $\BU_\alpha$ denote the one-dimensional unipotent subgroup of $\BG$
corresponding to $\alpha$ and let $x_\alpha : \bar{\BF}_p \to \BU_\alpha$ be an
isomorphism
of algebraic groups. We may, and we will, choose the family
$(x_\alpha)_{\alpha \in \Phi}$
such that $x_\alpha(\xi^{q_\alpha^\circ})=F(x_\alpha(\xi))$
for all $\alpha \in \Phi$ and $\xi \in \bar{\BF}_p$. In particular,
$F^{d_\alpha}(x_\alpha(\xi))=x_\alpha(\xi^{q_\alpha})$ and we deduce an
isomorphism $\BU_\alpha^{F^{d_\alpha}} \iso \BF_{q_\alpha}$.

\smallskip
We put $G=\BG^F$, $T_0=\BT_0^F$, and $U_0=\BU_0^F$.

\subsubsection{}
Let $I \subset \Delta$. We denote by
\begin{itemize}
\item $W_I$ the subgroup of $W$ generated by $(s_\alpha)_{\alpha \in I}$,
\item $W^I$ the set of elements $w \in W$ which are of minimal length
in $W_I w$,
\item $w_I$ the longest element of $W_I$,
\item $\BP_I$ the parabolic subgroup $\BB_0 W_I\BB_0$ of $\BG$,
\item $\BL_I$ the unique Levi complement of $\BP_I$ containing $\BT_0$,
\item $\BV_I$ the unipotent radical of $\BP_I$, 
\item $\BB_I$ the Borel subgroup $\BB \cap \BL_I$ of $\BL_I$,
\item $\BU_I$ the unipotent radical of $\BB_I$.
\end{itemize}

In particular, $\BU=\BU_I \ltimes \BV_I$.
If $I$ is $\phi$-stable, then $W_I$, $\BP_I$, $\BL_I$, $\BV_I$, $\BB_I$, and
$\BU_I$
are $F$-stable. Note that $W_\Delta=W$, $\BP_\Delta=\BL_\Delta=\BG$,
$\BB_\Delta=\BB_0$, $\BU_\Delta=\BU_0$ and $\BV_\Delta=1$.
On the other hand, $W_\emptyset=1$, $\BP_\emptyset=\BB_0$,
$\BL_\emptyset=\BB_\emptyset=\BT_0$, $\BU_\emptyset=1$, and
$\BV_\emptyset=\BU_0$.

We put $L_I=\BL_I^F$, etc...

\subsubsection{}
Let $D(\BU_0)$ be the derived subgroup of $\BU_0$. For any total
order on $\Phi_+$,
the product map $\prod_{\alpha \in \Phi^+ \setminus \Delta}
\BU_\alpha \to D(\BU_0)$ is an isomorphism
of varieties. It is not in general an isomorphism of algebraic groups.
The canonical map $\prod_{\alpha \in \Delta} \BU_\alpha \to \BU_0/D(\BU_0)$ is
an isomorphism of algebraic groups commuting with $F$.
We deduce an isomorphism from 
$(\BU_0/D(\BU_0))^F$ (which is canonically isomorphic to
$\BU_0^F/D(\BU_0)^F$) to
$\prod_{\alpha \in [\Delta/\phi]} \BU_\alpha^{F^{d_\alpha}}$,
which identifies with $\prod_{\alpha \in [\Delta/\phi]} \BF_{q_\alpha}$.

%We assume there is $f\in\{1,2\}$
%such that $F^f$ is a Frobenius endomorphism of
%$\BG$ over some finite field $\BF_{q_f}$.
%We denote by $q$ the positive real number with $q^f=q_f$.

\subsection{Deligne-Lusztig induction and restriction}
\subsubsection{Harish-Chandra induction and restriction}
Let $\BP$ be an $F$-stable parabolic subgroup of $\BG$, let $\BL$
be an $F$-stable Levi complement of $\BP$ and let $\BV$ denote the
unipotent radical of $\BP$. We put $L=\BL^F$, $P=\BP^F$, and $V=\BV^F$.
The {\it Harish-Chandra restriction} is the functor
\begin{align*}
{^*\CR}_{\BL\subset\BP}^{\BG}: \CO G\mMod&\to \CO L\mMod\\
M&\mapsto M^V.
\end{align*}
The {\it Harish-Chandra induction} is the functor
\begin{align*}
\CR_{\BL\subset\BP}^\BG: \CO L\mMod&\to \CO G\mMod\\
M&\mapsto\Ind_P^G \circ \Res_P^L M,
\end{align*}
where $\Res_P^L$ is defined through the canonical surjective
morphism $P \to L$. These functors are left and right adjoint
to each other (note that $|V|$ is invertible in $\CO$).

\subsubsection{Definition}
\label{secdefDL}
Let $\BP$ be a parabolic subgroup of $G$, let $\BL$ be a Levi
complement of $\BP$ and let $\BV$ denote the unipotent radical of $\BP$.
We assume that $\BL$ is $F$-stable.
Let
\begin{eqnarray*}
Y_\BV\!\!\!&=&\!\!\!Y_{\BV,\BG}=\{g\BV \in \BG/\BV~|~g^{-1}F(g) \in \BV\cdot F(\BV)\}\\
\textrm{and}\quad X_\BV\!\!\!&=&\!\!\!X_{\BV,\BG}=
\{g\BP \in \BG/\BP~|~g^{-1}F(g) \in \BP\cdot F(\BP)\}.
\end{eqnarray*}
Let $\pi_\BV : Y_\BV \to X_\BV$, $g\BV \mapsto g\BP$ be the canonical map.
The group $G$ acts on $Y_\BV$ and $X_\BV$
by left multiplication and $L$ acts on $Y_\BV$ by right
multiplication. Moreover, $\pi_\BV$ is $G$-equivariant and it is
the quotient morphism by $L$. Note that $Y_\BV$ and $X_\BV$ are smooth.

Let us consider $R\Gamma_c(Y_\BV,\CO)$, an element of
$D^b((\CO G)\otimes (\CO L)^\oppose)$, which is perfect for $\CO G$
and for $(\CO L)^\oppose$.
We have associated left and right adjoint functors
between the derived categories $D^b(\CO L)$ and $D^b(\CO G)$:
\begin{eqnarray*}
\CR_{\BL \subset \BP}^\BG\!\!\!&=&\!\!\!R\Gamma_c(Y_\BV,\CO)\otimes_{\CO L}^\BL-:
D^b(\CO L)\to D^b(\CO G)\\
\textrm{and}\quad 
{^*\CR}_{\BL \subset \BP}^\BG\!\!\!&=&\!\!\!
R\Hom^\bullet_{\CO G}(R\Gamma_c(Y_\BV,\CO),-):
D^b(\CO G)\to D^b(\CO L).
\end{eqnarray*}
Tensoring by $K$, they induce adjoint morphisms between Grothendieck groups:
$$R_{\BL \subset \BP}^\BG : \mathrm{K}_0(K L) \to \mathrm{K}_0(K G)
\textrm{ and }{^*R}_{\BL \subset \BP}^\BG :
\mathrm{K}_0(K G) \to \mathrm{K}_0(K L).$$

Note that when $\BP$ is $F$-stable, then the Deligne-Lusztig
functors are induced by the corresponding Harish-Chandra functors between
module categories.

\subsubsection{Reductions}
\label{section reductions}
We describe here some relations between 
Deligne-Lusztig varieties of groups of the same type. 

\smallskip
Let $\tilde{\BG}$ be a connected reductive algebraic group over
$\bar{\BF}_p$ endowed with an endomorphism $F$ and assume a
non-trivial power of $F$ is a Frobenius endomorphism.
Let $i:\BG\to\tilde{\BG}$ be a morphism of algebraic groups commuting with
$F$ and such that the kernel $Z$ of $i$
is a central subgroup of $\BG$ and the image of $i$ contains
$[\tilde{\BG},\tilde{\BG}]$.

Let $\tilde{\BP}$ be a parabolic subgroup of $\tilde{\BG}$ with an
$\tilde{F}$-stable Levi complement $\tilde{\BL}$. Let
$\BP=i^{-1}(\tilde{\BP})$ and $\BL=i^{-1}(\tilde{\BL})$, a parabolic
subgroup of $\BG$ with an $F$-stable Levi complement.
Let $\tilde{\BV}$ (resp. $\BV$)
be the unipotent radical of $\tilde{\BP}$ (resp. $\BP$).

Let $N=\{(x,l)\in \tilde{\BG}^F\times (\tilde{\BL}^F)^\oppose|
xl^{-1}\in i(\BG^F)\}$
and consider its action on $\tilde{\BG}^F$ given by $(x,l)\cdot g=xgl^{-1}$.
Then, $N$ stabilizes $i(\BG^F)$.
The morphism $i$ induces a canonical morphism
$Y_{\BV,\BG}\to Y_{\tilde{\BV},\tilde{\BG}}$. Its image is isomorphic to
$(Z^F\times (Z^F)^\oppose)\backslash Y_{\BV,\BG}=
Y_{\BV,\BG}/Z^F=Z^F\backslash Y_{\BV,\BG}$ and it is stable under the
action of $N$.

\begin{prop}
\label{isobetweenY}
The canonical map induces a radicial
$(\tilde{\BG}^F\times (\tilde{\BL}^F)^\oppose)$-equivariant map
$$\Ind_N^{\tilde{\BG}^F\times (\tilde{\BL}^F)^\oppose}(Y_{\BV,\BG}/Z^F)\to
Y_{\tilde{\BV},\tilde{\BG}},$$
\ie, we have radicial maps
$$\tilde{\BG}^F\times_{i(\BG^F)}Y_{\BV,\BG}\leftarrow
Y_{\tilde{\BV},\tilde{\BG}}\to
Y_{\BV,\BG}\times_{i(\BL^F)}\tilde{\BL}^F$$
and we deduce a canonical isomorphism in
$K^b(\CO(\tilde{\BG}^F\times (\tilde{\BL}^F)^\oppose))$
$$\Ind_N^{\tilde{\BG}^F\times (\tilde{\BL}^F)^\oppose}
\tilde{R}\Gamma_c(Y_{\BV,\BG},\CO)^{Z^F}\iso
\tilde{R}\Gamma_c(Y_{\tilde{\BV},\tilde{\BG}},\CO).$$
%We have isomorphisms of functors
%$$i_* \CR_{\BL\subset\CP}^\BG \osi
% \CR_{\tilde{\BL}\subset\tilde{\BP}}^{\tilde{\BG}}\iso
%\CR_{\BL\subset\CP}^\BG$$
\end{prop}

\begin{proof}
Since $Y_{\tilde{\BV},\tilde{\BG}}$ is smooth, it is enough to prove that
the map is bijective on (closed) points to deduce that it is radicial
(cf e.g. \cite[AG.18.2]{Bor}).

We can factor $i$ as the composition of three maps:
$\BG\to \BG/Z^0 \to \BG/Z\to \tilde{\BG}$.
The finite group $Z/Z^0$ has an $F$-invariant
filtration $1=L^0\subset L^1\subset\cdots L^n=Z/Z^0$ where
$F$ acts trivially on $L^i/L^{i+1}$ for $0\le i<n$ and
the Lang map is an isomorphism on $L^n/L^{n-1}$. This allows us
to reduce the proof of the Proposition to one of the following three cases:
\begin{itemize}
\item[(1)] $i$ is surjective and the Lang map on $Z$ is surjective
\item[(2)] $i$ is surjective and $Z^F=Z$
\item[(3)] $i$ is injective.
\end{itemize}
The first and third cases are solved as in \cite[Lemma 2.1.2]{Bon}. Let us
consider the second case (cf \cite[Proof of Proposition 2.2.2]{Bon}).
Let $\tau:\tilde{\BG}^F/i(\BG^F)\iso Z$ be the canonical isomorphism: given
$g\in\tilde{\BG}^F$, let $h\in \BG$ with $i(h)=g$. Then,
$\tau(g)=h^{-1}F(h)$.

We need to show that $Y_{\tilde{\BV},\tilde{\BG}}=
\coprod_{f\in \tilde{\BG}^F/i(\BG^F)}f\cdot i(Y_{\BV,\BG})$.
Let $y\in Y_{\tilde{\BV},\tilde{\BG}}$. Let $g\in \BG$ such that $y=i(g)$.
There is $z\in Z$ such that $g^{-1}F(g)\in z\BV\cdot F(\BV)$.
Let now $f\in\tilde{\BG}^F$
with $\tau(f)=z^{-1}$ and $g'\in \BG$ with $f=i(g')$. We have
$g'g\in Y_{\BV,\BG}$, hence $fy\in i(Y_{\BV,\BG})$ and we are done.

\smallskip
The assertion about $\tilde{R}\Gamma_c$ follows from the fact that this
complex depends
only on the \'etale site and a radicial morphism induces an equivalence of 
\'etale
sites and from the fact that
$\tilde{R}\Gamma_c$ of a quotient is isomorphic to the fixed points on
$\tilde{R}\Gamma_c$ (cf \cite[Theorem 4.1]{Ri1} or
\cite[Th\'eor\`eme 2.28]{Rou2}).
\end{proof}

Note that $i$ induces an isomorphism
$X_{\BV,\BG}\iso X_{\tilde{\BV},\tilde{\BG}}$ and the isomorphisms of
Proposition \ref{isobetweenY} are compatible with that isomorphism.

\subsubsection{Local study}
\label{seclocal}
We keep the notation of \S \ref{secdefDL}. 
Recall (\cite{Rou2}, \cite{Ri1}) that $\tilde{R}\Gamma_c(Y_\BV,\CO)$ is
a bounded complex of
$\CO(G\times L^\oppose)$-modules that
is homotopy equivalent to a bounded complex
$C'$ of finitely generated modules which are direct summands of permutation
modules. We can furthermore assume that $C'$ has no direct summand homotopy
equivalent to $0$. Let $S$ be the Sylow $\ell$-subgroup of $Z(L)$.
Let $b$ be the sum of the block idempotents of $\CO G$ whose defect group
contains (up to conjugacy) $S$.

\begin{lemma}
\label{local}
If $C_\BG(S)=\BL$, then the canonical morphism $b\CO G\to
\End_{K^b(\CO L^\oppose)}(bC')$
is a split injection of $(b\CO G\otimes (b\CO G)^\oppose)$-modules.
\end{lemma}

\begin{proof}
We have a canonical isomorphism $\End_{\CO L^\oppose}^\bullet(bC')\iso
bC'\otimes_{\CO L}\Hom_{\CO}^\bullet(bC',\CO)$ and a canonical map
$C'\otimes_{\CO L}\Hom_{\CO}^\bullet(bC',\CO)\to b\CO G$. We will show that the
composition
$$b\CO G\to \End_{\CO L^\oppose}^\bullet(bC')\to b\CO G$$
is an isomorphism, which will suffice, since
$\End_{K^b(\CO L^\oppose)}(bC')=H^0(\End_{\CO L^\oppose}^\bullet(bC'))$.

Let $e$ be a block idempotent of $b\CO G$. We will show (and this will suffice) 
that the composition above is an isomorphism after multiplication by $e$ and
after tensoring by $k$
$$ekG\xrightarrow{\alpha} \End_{kL^\oppose}^\bullet(ekC')
\xrightarrow{\beta} ekG.$$

The composition $\beta\alpha$ is an endomorphism of $(ekG,ekG)$-bimodules of
$ekG$, \ie, an element of the local algebra $Z(ekG)$. So, it suffices to
show that $\beta\alpha$ is not nilpotent. This will follow from the
non-nilpotence of its image under $\Br_{\Delta S}$, where
$\Delta S=\{(x,x^{-1})|x\in S\}\subset G\times L^\oppose$.

By \cite[Proofs of Lemma 4.2 and Theorem 4.1]{Ri2}, there is a commutative
diagram
$$\xymatrix{
\br_S(e)kL \ar[rr]^-{\Br_{\Delta S}(\alpha)}\ar[drr]_{a} &&
 \Br_{\Delta S}(\End_{kL^\oppose}^\bullet(eC'))\ar[rr]^-{\Br_{\Delta S}(\beta)} &&
\br_S(e)kL \\
&& \End_{kL^\oppose}^\bullet(\Br_{\Delta S}(eC'))\ar[urr]_{b}\ar[u]
}$$
where $a$ and $b$ are the natural maps coming from the left action of
$\br_S(e)kL$ on $\Br_{\Delta S}(eC')=\br_S(e)\Br_{\Delta S}(C')$. Note
that $\br_S(e)\not=0$.

We have $\Br_{\Delta S}(C')\simeq kL$ (cf \cite[proof of Theorem 4.5]{Rou2},
the key point is that $Y_\BV^{\Delta S}=L$).
So, $a$ and $b$ are isomorphisms.
\end{proof}

\subsection{Alternative construction for tori}
It will be convenient to use an alternative description of Deligne-Lusztig
functors in the case of tori (as in \cite[\S 11.2]{BonRou} for example).

\subsubsection{}
\label{modelsRTG}
Let $w \in W$. We define
\begin{eqnarray*}
Y(\dot{w})\!\!\!&=&\!\!\!Y_{\BG,F}(\dot{w}) =
 \{g\BU_0 \in \BG/\BU_0~|~g^{-1}F(g) \in \BU_0 \dot{w} \BU_0 \} \\
\textrm{and}\quad 
X(w)\!\!\!&=&\!\!\!X_{\BG,F}(w)
= \{g\BB_0 \in \BG/\BB_0~|~g^{-1}F(g) \in \BB_0 w\BB_0 \}
\end{eqnarray*}
Let $\pi_w : Y(\dot{w}) \to X(w)$, $g\BU_0 \mapsto g\BB_0$. The group $G$
acts on
$Y(\dot{w})$ and $X(w)$ by left multiplication while $\BT_0^{wF}$ acts on
$Y(\dot{w})$
by right multiplication. Moreover, $\pi_w$ is $G$-equivariant and it is
isomorphic to the quotient morphism by $\BT_0^{wF}$.
The complex $R\Gamma_c(Y(\dot{w}),\CO)$, viewed in
$D^b((\CO G)\otimes (\CO \BT_0^{wF})^\oppose)$, induces
left and right adjoint functors
between the derived categories $D^b(\CO \BT_0^{wF})$ and $D^b(\CO G)$:
\begin{eqnarray*}
\CR_{\dot{w}}\!\!&=&\!\!R\Gamma_c(Y(\dot{w}),\CO)\otimes_{\CO \BT_0^{wF}}^\BL-:
D^b(\CO \BT_0^{wF})\to D^b(\CO G)\\
\textrm{and}\quad
{^*\CR}_{\dot{w}}\!\!&=&\!\!R\Hom^\bullet_{\CO G}(R\Gamma_c(Y(\dot{w}),\CO),-):
D^b(\CO G)\to D^b(\CO \BT_0^{wF}).
\end{eqnarray*}

\medskip
Let $g \in \BG$ be such that $g^{-1} F(g)=\dot{w}$.
Let $\BT_w=g\BT_0 g^{-1}$, $\BB_w=g \BB_0 g^{-1}$ and $\BU_w = g\BU_0g^{-1}$.
Then, $\BT_w$ is an $F$-stable maximal torus of the Borel subgroup $\BB_w$ of
$\BG$, conjugation by $g$ induces an isomorphism
$\BT_0^{wF} \simeq \BT_w^F$, and
multiplication by $g^{-1}$ on the right induces an isomorphism
$Y(\dot{w}) \iso Y_{\BU_w}$
which is equivariant for the actions of $G$ and $\BT_0^{wF}\simeq \BT_w^F$.
This provides an identification between the Deligne-Lusztig induction functors
$\CR_{\BT_w \subset \BB_w}^\BG$ and $\CR_{\dot{w}}$
(and similarly for the Deligne-Lusztig restriction functors).

\subsubsection{}
Let us consider here the case of a product of groups cyclically permuted by
$F$ (cf \cite[\S 1.18]{Lu2}, \cite[Proposition 2.3.3]{DiMiRou}, whose
proofs for the varieties $X$ extend to the varieties $Y$).

\begin{prop}
\label{cycles}
Assume $\BG=\BG_1\times\cdots\times\BG_s$ with $F(\BG_i)=\BG_{i+1}$ for
$i<s$ and $F(\BG_s)=\BG_1$.
Let $\BB_1$ be an $F^s$-stable Borel subgroup of $\BG_1$ and
$\BT_1$ an $F^s$-stable maximal torus of $\BB_1$. Let
$\BB_i=F^i(\BB_1)$ and $\BT_i=F^i(\BT_1)$. We assume that
$\BB_0=\BB_1\times\cdots\times\BB_s$ and $\BT_0=\BT_1\times\cdots\times\BT_s$.
Let $W_i$ be the Weyl group of $\BG_i$ relative to $\BT_i$. We
identify $W$ with $W_1\times\cdots\times W_s$.
Let $w=(w_1,w_2,\ldots,w_s)\in W$ and let
$v=F^{1-s}(w_s)F^{2-s}(w_{s-1})\cdots w_1$. 
The first projection gives isomorphisms
$\BG^F\iso\BG_1^{F^s}$ and
$\BT_0^{wF}\iso \BT_1^{vF^s}$ and we identify
the groups via these isomorphisms.

If $l(w_1)+\cdots+l(w_s)=l(v)$, then there is a canonical 
$(\BG^F\times (\BT_0^{wF})^\oppose)$-equivariant radicial map
$$Y_{\BG,F}(\dot{w})\iso Y_{\BG_1,F^s}(\dot{w}_1).$$
\end{prop}

%\begin{rem}
%\label{remorbit}
%For a general element $w=(w_1,\ldots,w_s)$ of $W$, the previous Proposition
%generalizes to an isomorphism
%between $Y_{\BG,F}(\dot{w})$ and a generalized Deligne-Lusztig variety 
%$Y_{\BG_1,F^s}()$.
%Let $v=F^{1-s}(w_s)F^{2-s}(w_{s-1})\dots w_1$
%If $l(w_1)+\cdots+l(w_s)=l(v)$, then we have a canonical isomorphism
%and we obtain an isomorphism
%$$Y_{\BG,F}(\dot{w})\iso Y_{\BG_1,F^s}(v).$$
%\end{rem}

\begin{rem}
In general, let $\hat{\BG}$ be the simply connected covering of $[\BG,\BG]$.
This provides a morphism $i:\hat{\BG}\to\BG$ for which Proposition
\ref{isobetweenY} applies. This reduces the study of the variety $Y(\dot{w})$
to the case of a simply connected group. Assume now $\BG$ is simply connected.
Then, there is a decomposition $\BG=\BG_1\times\cdots\times\BG_s$ where each
$\BG_i$ is quasi-simple and $F$ permutes the components.
A variety $Y(\dot{w})$ for $\BG$ decomposes as a product of varieties for
the products
of groups in each orbit of $F$. So, the study of $Y(\dot{w})$ is further
reduced to the case where $F$ permutes transitively the components.
Proposition \ref{cycles} reduces finally the study to the case where
$\BG$ is quasi-simple, when $w$ is of the special type described in Proposition 
\ref{cycles}.
\end{rem}

\subsubsection{}

\begin{prop}
\label{Yestconnexe}
Let $w \in W$.
Assume that $\BG$ is semisimple and simply-connected. Then,
$Y(\dot{w})$ is irreducible if and only if $X(w)$ is irreducible.
\end{prop}

\begin{proof}
First, since $\pi_w : Y(\dot{w}) \to X(w)$ is surjective, 
we see that if $Y(\dot{w})$ is irreducible, then $X(w)$ is irreducible.

Let us assume now that $X(w)$ is irreducible.
Since $X(w)=Y(\dot{w})/\BT_0^{wF}$ is irreducible, it follows that $\BT_0^{wF}$ 
permutes transitively the irreducible components of $Y(\dot{w})$. Note that 
$\BG^F\backslash Y(\dot{w})$ is also irreducible, hence $\BG^F$ permutes
transitively the irreducible components of $Y(\dot{w})$. Let $Y^\circ$ be
an irreducible component of $Y(\dot{w})$, 
$H_1$ its stabilizer in $\BT_0^{wF}$ and $H_2$ its stabilizer in $\BG^F$.
Since the actions of $\BT_0^{wF}$ and $\BG^F$ on $Y(\dot{w})$ commute, it follows
that $H_2$ is a normal subgroup of $\BG^F$ and $\BG^F/H_2\simeq \BT_0^{wF}/H_1$.
But $G^F/[G^F,G^F]$ is a $p$-group because $\BG$ is semi-simple and simply 
connected \cite[Theorem 12.4]{St}. On the other hand, 
$\BT_0^{wF}/H_1$ is a $p'$-group. So $H_2=\BG^F$ and we are done.
\end{proof}

\begin{rem}
One can actually show that 
$X(w)$ is irreducible if and only if $w$ is not contained 
in a proper standard $F$-stable parabolic subgroup of $W$ (see \cite{cras}).
\end{rem}

\section{Gelfand-Graev representations}
\subsection{Definitions}
\subsubsection{}
Let $\psi : U_0 \to \Lambda^\times$ be a linear character trivial on
$D(\BU_0)^F$ and let $\alpha \in [\Delta/\phi]$. We denote by
$\psi_\alpha : \BF_{q_\alpha} \to \Lambda^\times$
the restriction of $\psi$ through
$\BF_{q_\alpha} \iso \BU_\alpha^{F^{d_\alpha}} \hookrightarrow
\BU_0^F/D(\BU_0)^F$.
Conversely, given
$(\psi_\alpha : \BF_{q_\alpha} \to \Lambda^\times)_{\alpha \in [\Delta/\phi]}$
a family of linear characters, we denote by
$\boxtimes_{\alpha \in [\Delta/\phi]} \psi_\alpha$
the linear character $U_0 \to \Lambda^\times$ trivial on $D(\BU_0)^F$
such that $(\boxtimes_{\alpha \in [\Delta/\phi]} \psi_\alpha)_\beta =
\psi_\beta$ for every $\beta \in [\Delta/\phi]$.

A linear character $\psi$ is {\it regular} (or $G$-regular if we need
to emphasize the ambient group) if it is trivial on $D(\BU_0)^F$ and if
$\psi_\alpha$ is non-trivial for every $\alpha \in [\Delta/\phi]$.

\subsubsection{Levi subgroups}
\label{psi I} 
Let $\psi$ be a linear character of $U_0$ which is trivial on $D(\BU_0)^F$. 
If $I$ is a $\phi$-stable subset of $\Delta$, we set $\psi_I=\Res_{U_I}^{U_0}\psi$. 
If $\psi'$ is a linear character of $U_I$ which is trivial on $D(\BU_I)^F$, we 
denote by $\tilde{\psi}'$ the restriction of $\psi'$ through the morphism 
$U_0 \to U_0/V_I\iso U_I$. This is a linear character of $U_0$ 
which is trivial on $D(\BU_0)^F$. 

\subsubsection{}
Since $U_0$ is a $p$-group, given $R=K$ or $R=k$, the canonical map
$\Hom(U_0,\CO^\times) \to \Hom(U_0,R^\times)$ is an isomorphism and
these groups will be identified.
Given $\psi \in \Hom(U_0,\CO^\times)$, we set
$$\Gamma_\psi=\Ind_{U_0}^G \CO_\psi.$$
If the ambient group is not clear from the context, $\Gamma_\psi$
will be denoted by $\Gamma_{\psi,G}$. 
Note that $\Gamma_\psi$ is a projective $\CO G$-module.
If $\psi$ is regular, then $\Gamma_\psi$ is a {\it Gelfand-Graev module}
of $G$ and the character of $K\Gamma_\psi$ is called a
{\it Gelfand-Graev character}.

\subsection{Induction and restriction}
\label{HC}
The results in this \S \ref{HC} are all classical.

\subsubsection{Harish-Chandra restriction}
The next Proposition
shows the Harish-Chandra restriction of a Gelfand-Graev
module is a Gelfand-Graev module, a result
of Rodier over $K$ (see \cite[Proposition 8.1.6]{Ca}). It must be noticed
that, in \cite[Chapter 8]{Ca}, the author works under the hypothesis that 
the centre is connected. However, the proof of Rodier's result 
given in \cite[Proof of Proposition 8.1.6]{Ca} remains valid when
the centre is not connected (see for instance \cite[Theorem 2.9]{DiLeMi2}). 
In both cases, the proof proposed is character-theoretic. 
Since we want to work with modular representations, we present 
here a module-theoretic argument, which is only the translation of the 
previous proofs.

\begin{thm}
\label{hc restriction}
Let $\psi : U_0 \to \CO^\times$ be a regular linear character.
Let $I$ be a $\phi$-stable subset of $\Delta$.
Let $\psi_I=\Res_{U_I}^{U_0} \psi$.
Then $\psi_I$ is an $L_I$-regular linear character and
$${^*\CR}_{\BL_I\subset\BP_I}^\BG \Gamma_{\psi,G} \simeq
\Gamma_{\psi_I,L_I}.$$
\end{thm}

\begin{proof}
As explained above, this result is known over
$K$ using scalar products of characters. The result is an
immediate consequence, since ${^*\CR}_{\BL_I}^\BG \Gamma_{\psi,G}$ is
projective and two projective modules with equal characters are isomorphic.
We provide nevertheless here a direct module-theoretic
proof --- it shows there is a canonical isomorphism.

First, it is clear from the definition that $\psi_I$ is a regular linear
character of $U_I$. Moreover,
$${^*\CR}_{\BL_I\subset\BP_I}^\BG \Gamma_{\psi,G} \simeq
\Bigl(\Res_{P_I}^{G} \Ind_{U_0}^{G} \CO_\psi \Bigr)^{V_I}.$$
The map
$(W^I)^F \to P_I \backslash G / U_0$,
$w \mapsto P_I w U_0$ is bijective. Thus,
by the Mackey formula for classical induction and restriction, we have
$${^*\CR}_{\BL_I\subset\BP_I}^\BG \Gamma_{\psi,G} \simeq
\bigoplus_{w \in (W^I)^F} \Bigl(\Ind_{P_I \cap {^wU}_0}^{P_I}
\Res_{P_I \cap {^wU}_0}^{{^wU}_0}
\CO_{^{\dot{w}}\psi}\Bigr)^{V_I}.$$
Fix $w \in (W^I)^F$. Then, $\BP_I \cap {^w\BU}_0 =
(\BL_I \cap {^w\BU}_0)\cdot(\BV_I \cap {^w\BU}_0)$ (see for instance
\cite[2.9.3]{DiLeMi1}). So, if the restriction of ${^w\psi}$ to
$V_I \cap {^wU}_0$ is non-trivial, then
$$\Bigl(\Ind_{P_I \cap {^wU}_0}^{P_I}
\Res_{P_I \cap {^wU}_0}^{{^wU}_0}
\CO_{^{\dot{w}}\psi}\Bigr)^{V_I}=0.$$
By \cite[Page 163]{DiLeMi1}, this happens unless $w = w_I w_\Delta$.
On the other hand, $\BP_I \cap {^{w_I w_\Delta}\BU}_0 = \BU_I$. Therefore,
$${^*\CR}_{\BL_I\subset\BP_I}^\BG \Gamma_{\psi,G} =
\Bigl(\Ind_{U_I}^{P_I} \Res_{U_I}^{{^{w_I w_\Delta}U}_0}
\CO_{^{\dot{w}_I\dot{w}_\Delta}\psi}\Bigr)^{V_I}.$$
In other words,
$${^*\CR}_{\BL_I\subset\BP_I}^\BG \Gamma_{\psi,G} \simeq
\Ind_{U_I}^{L_I} \Res_{U_I}^{{^{w_I w_\Delta}U}_0}
\CO_{^{\dot{w}_I\dot{w}_\Delta}\psi}.$$
Since $\psi_I$ is $T_0$-conjugate to $\Res_{U_I}^{^{w_I w_\Delta}{U}_0}
({^{\dot{w}_I\dot{w}_\Delta}\psi})$ by \cite[2.9.6]{DiLeMi1}, we get
the result.
\end{proof}

\subsubsection{Harish-Chandra induction}
\label{hc induction}

Let $I$ be a $\phi$-stable subset of $\Delta$ and let
$\psi : U_I \to \CO^\times$ be a linear character of $U_I$. Then
$$\CR_{\BL_I\subset\BP_I}^\BG \Gamma_{\psi,L_I} \simeq
\Gamma_{\tilde{\psi},G},$$
where $\tilde{\psi} : U_0 \to \CO^\times$ is the lift of $\psi$ through the 
canonical map $U_0 \to U_0/V_I \iso U_I$.

\subsection{Deligne-Lusztig restriction}
\label{secconjectures}
\subsubsection{Restriction to Levi subgroups}
Let $\BP$ be a parabolic subgroup of $\BG$ with an $F$-stable Levi
complement $\BL$. Let $\BV$ be the unipotent radical of $\BP$. 
We denote by $\BU_\BL$ the unipotent radical of some $F$-stable
Borel subgroup of $\BL$. The proof of the next result is due 
to Digne, Lehrer and Michel. If the centre of $\BL$ is disconnected, 
then the proof is given in \cite{DiLeMi2}: it requires the theory of 
character sheaves. This explains why the scope of validity of this 
result is not complete, but it is reasonable to hope that it holds 
in general. If the centre of $\BL$ is connected, then see \cite{DiLeMi1}.

\begin{thm}[Digne-Lehrer-Michel]
\label{dl restriction}
Assume that one of the following holds:
\begin{itemize}
\item[(1)] $p$ is good for $G$, $q$ is large enough and $F$ is a Frobenius 
endomorphism;

\item[(2)] the center of $\BL$ is connected.
\end{itemize}
Let $\psi : U_0 \to \CO^\times$ be a $G$-regular linear character. Then there
exists an $L$-regular linear character $\psi_L$ of $U_L$ such that
$${^*R}_{\BL \subset \BP}^\BG [K\Gamma_{\psi,G}]=(-1)^{\dim Y_\BV} 
[K\Gamma_{\psi_L,L}].$$
\end{thm}

\bigskip

We have good evidences that a much stronger result holds:

\begin{conj}
\label{conj1}
 Let $\psi : U_0 \to \CO^\times$ be a $G$-regular
linear character. Then there exists an $L$-regular linear character
$\psi_L$ of $U_L$ such that
$${^*\CR}_{\BL \subset \BP}^{\BG} \Gamma_{\psi,G}\simeq
\Gamma_{\psi_L,L}[-\dim Y_\BV].$$
\end{conj}

It is immediate that Conjecture \ref{conj1}
is compatible with Theorem \ref{dl restriction}. If $\BP$ is $F$-stable, then
the conjecture holds by Theorem \ref{hc restriction}.
As we will see in Theorem \ref{coxeter}, the conjecture holds if $\BL$ is a
maximal torus and $(\BP,F(\BP))$ lies in the $\BG$-orbit of
$(\BB_0,{^w\BB}_0)$, where $w$ is a product of simple reflections lying in
different $F$-orbits.

\smallskip
Note that Conjecture \ref{conj1}
is also compatible with the Jordan decomposition
\cite[Theorem B']{BonRou}.

\begin{rem}
Let us examine the consequences of Conjecture \ref{conj1} at the level
of $KG$-modules. An irreducible
character of $G$ is called {\it regular} if it is a component of a
Gelfand-Graev character of $G$ (for instance, the Steinberg character
is regular). Then, Conjecture \ref{conj1} over $K$ is equivalent to the
statement that regular characters
appear only in degree $\dim Y_\BV$ in the cohomology of the Deligne-Lusztig
variety $Y_\BV$. This is known for the Steinberg character,
if $\BL$ is a maximal torus \cite[Proposition 3.3.15]{DiMiRou}.
\end{rem}

\subsubsection{Restriction to tori}
We now restate Conjecture \ref{conj1} in the case of tori:

\begin{conj}
\label{conj2} Let $w \in W$ and let $\psi : U_0 \to \CO^\times$ be a regular
linear character. Then
$${^*\CR}_{\dot{w}} \Gamma_{\psi,G}\simeq\CO \BT_0^{wF}[-l(w)].$$
\end{conj}

\begin{rem}
Through the identification of \S \ref{modelsRTG}, Conjecture \ref{conj2} 
is equivalent to Conjecture \ref{conj1} for tori. Indeed, 
$\CO \BT_w^F \simeq \CO \BT_0^{wF}$ is the unique Gelfand-Graev module
of $\BT_w^F$.
\end{rem}

\bigskip

We now propose a conjecture which refines Conjecture \ref{conj2}.
We first need some notation. Given $x,w \in W$, we put
\begin{eqnarray*}
Y_x(\dot{w})\!\! &=&\!\! \{g\BU_0 \in \BB_0x\BU_0/\BU_0~|~g^{-1}F(g) \in
 \BU_0 \dot{w} \BU_0 \}\\
\textrm{and}\quad X_x(w)\!\! &=&\!\! \{g\BB_0 \in \BB_0x\BB_0/\BB_0~|~g^{-1}F(g) \in
 \BB_0 w\BB_0 \}.
\end{eqnarray*}
Then $(Y_x(\dot{w}))_{x \in W}$ (resp. $(X_x(w))_{x \in W}$)
is a stratification of $Y(\dot{w})$ (resp. $X(w)$). Note that some strata
might be empty. Note also that $Y_x(\dot{w})$ and $X_x(w)$ are stable under
left multiplication by $B_0$.
Moreover, $Y_x(\dot{w})$ is stable by right multiplication by $\BT_0^{wF}$ and 
the natural map $Y_x(\dot{w}) \to X_x(w)$, $g \BU_0 \to g\BB_0$ is the 
quotient morphism for this action.

\begin{conj}
\label{conj2plus}
Let $\psi : U_0 \to \CO^\times$ be a $G$-regular linear character. Then, there
are isomorphisms in $D^b(\CO \BT_0^{wF})$:
$$R\Hom^\bullet_{\CO U_0}(R\Gamma_c(Y_x(\dot{w}),\CO),\CO_\psi) \simeq
\begin{cases}
\CO \BT_0^{wF} & \text{if }x=w_\Delta,\\
0 & \text{otherwise.}
\end{cases}$$
\end{conj}

Note that 
${^*\CR}_{\dot{w}} \Gamma_{\psi,G}\simeq
R\Hom^\bullet_{\CO U_0}(R\Gamma_c(Y(\dot{w}),\CO),\CO_\psi)$, hence
Conjecture \ref{conj2plus} implies Conjecture \ref{conj2}.

\begin{rem}
\label{reduceLevi}
The proof of Theorem \ref{hc restriction} shows that that
Conjecture \ref{conj2plus} holds for $w=1$. More generally,
if $I$ is a $\phi$-stable subset of $\Delta$ and if $w \in W_I$, then
Conjecture \ref{conj2plus} holds for $(\BG,w)$ if and only if 
it holds for $(\BL_I,w)$.
\end{rem}

\section{Coxeter orbits}\label{section coxeter}

\begin{quotation}
\noindent{\bf Notation:} {\it Let $r=|\Delta/\phi|$. 
We write $[\Delta/\phi]=\{\alpha_1,\dots,\alpha_r\}$.
Let $w = s_{\alpha_1} \dots s_{\alpha_r}$, a {\it twisted Coxeter element}.
We put $T=\BT_0^{wF}$.}
\end{quotation}
The aim of this section is to study the complex of
cohomology $C=R\Gamma_c(Y(\dot{w}),\CO)$ of the Deligne-Lusztig 
variety $Y(\dot{w})$. As a consequence, we get that Conjecture 
\ref{conj2plus} holds for $w$ (see Theorem \ref{coxeter}). 

\subsection{The variety $D(\BU_0)^F\backslash X(w)$}
\subsubsection{}
The variety $X(w)$ has been studied by Lusztig \cite[\S 2]{Lu2}.
Before summarizing some of his results, we need some notation.
Let
$$X'(w)=X'_\BG(w)=
\{u \in \BU_0~|~u^{-1} F(u) \in \BU_{-w_\Delta(\alpha_1)}^\# \cdots
 \BU_{-w_\Delta(\alpha_r)}^\#\}.$$
Then, we have \cite[Corollary 2.5, Theorem 2.6 and Proposition 4.8]{Lu2}:

%\begin{thm}[Lusztig]
%\label{Lusztig}
%\begin{equation}\label{inclusion}
%X(w) \subset \BB_0 w_\Delta \BB_0/\BB_0.
%\end{equation}
%\begin{equation}
%\label{isomorphisme}
%{\text{\it $X'(w) \to X(w)$, $u \mapsto u w_\Delta \BB_0$ is an
%isomorphism of varieties.}}
%\end{equation}
%\begin{equation}
%\label{irreductible}
%{\text{\it The variety $X(w)$ is irreducible.}}
%\end{equation}
%\end{thm}

\begin{thm}[Lusztig]
\label{Lusztig}
We have:
\begin{itemize}
\item[$({\mathrm{a}})$] 
$X(w) \subset \BB_0 w_\Delta \BB_0/\BB_0$. 

\item[$({\mathrm{b}})$] The map $X'(w) \to X(w)$, $u \mapsto u w_\Delta \BB_0$ is an
isomorphism of varieties.

\item[$({\mathrm{c}})$] The variety $X(w)$ is irreducible.
\end{itemize}
\end{thm}

The isomorphism (b) above is $B_0$-equivariant ($U_0$ acts on
$\BU_0$ by left multiplication and $T_0$ acts on $\BU_0$ by conjugation).

The Lang map $\BU_0 \to \BU_0$, $u \mapsto u^{-1}F(u)$, induces an
isomorphism
$$U_0 \backslash X'(w) \iso \BU_{-w_\Delta(\alpha_1)}^\# \cdots
 \BU_{-w_\Delta(\alpha_r)}^\# \iso (\BG_m)^r.$$
Let
$$X''(w)=\{u D(\BU_0) \in \BU_0 / D(\BU_0)~|~u^{-1} F(u) \in
\BU_{-w_\Delta(\alpha_1)}^\# \cdots \BU_{-w_\Delta(\alpha_r)}^\# D(\BU_0)\}.$$
Then, $X''(w)$ is an open subvariety of $\BU_0/D(\BU_0)$. Moreover,
the canonical map $X'(w) \to X''(w)$, $u \mapsto uD(\BU_0)$ factorizes
through $D(\BU_0)^F \backslash X'(w)$. In fact:

\begin{prop}\label{quotient par DU}
The induced map $D(\BU_0)^F \backslash X'(w) \to X''(w)$ is a $U_0$-equivariant
isomorphism of varieties.
\end{prop}

\begin{proof}
Let $Z=\BU_{-w_\Delta(\alpha_1)}^\# \cdots \BU_{-w_\Delta(\alpha_r)}^\#$ and
consider the morphism
$\tau : X'(w) \to X''(w)$, $u \mapsto u D(\BU_0)$. First, let us show that
the fibers of $\tau$ are $D(\BU_0)^F$-orbits. Let $u$ and $u'$ be two elements
of $X'(w)$ such that $\tau(u)=\tau(u')$. There exists $v \in D(\BU_0)$ such that
$u'=vu$. Therefore,
$$u^{\prime -1} F(u') = u^{-1} (v^{-1} F(v)) u \cdot u^{-1} F(u).$$
Recall that $u^{\prime -1} F(u')$ and $u^{-1} F(u)$ belong to
$Z$ and the map
$Z \to \BU_0/D(\BU_0)$ is injective.
Therefore, $u^{-1} (v^{-1} F(v)) u =1$. In other words, $F(v)=v$.

Let us now show that $\tau$ is surjective. Let $u \in \BU_0$ with
$u^{-1} F(u) \in Z D(\BU_0)$.
Write $u^{-1} F(u) = yx$, with $x \in Z$
and $y \in D(\BU_0)$. By Lang's Theorem applied to the group $D(\BU_0)$
and the isogeny $D(\BU_0) \to D(\BU_0)$, $v \mapsto x F(v) x^{-1}$, there
exists $v \in D(\BU_0)$ such that $v x F(v)^{-1} x^{-1} = y$. Then,
$(uv)^{-1} F(uv) = x$. So $uv \in X'(w)$ and $\tau(uv)=u D(\BU_0)$.
So $\tau$ is surjective.

Therefore, is it sufficient to show that $\tau$ is \'etale. Since
the maps $X'(w) \to U_0 \backslash X'(w)$ and
$X''(w) \to U_0 \backslash X''(w)=(\BU_0^F/D(\BU_0)^F) \backslash X''(w)$
are \'etale, it is sufficient to show that the map
$\tau' : U_0 \backslash X'(w) \to U_0 \backslash X''(w)$ induced by
$\tau$ is an isomorphism.
Via the canonical isomorphisms $U_0 \backslash X'(w) \iso Z$ and
$U_0 \backslash X''(w) \iso Z D(\BU_0) / D(\BU_0)$, then
$\tau'$ is the canonical map $Z\to Z D(\BU_0) / D(\BU_0)$,
which is clearly an isomorphism.
\end{proof}

\subsubsection{}
Let $i\in\{1,\ldots,r\}$, let $w_i=s_{\alpha_1}\cdots s_{\alpha_{i-1}}
s_{\alpha_{i+1}} \cdots s_{\alpha_r}$, and
let $I$ be the complement of the $\phi$-orbit of $\alpha_i$ in $\Delta$.
Let 
$$X'_i(w)=\{u \in \BU_0~|~u^{-1} F(u) \in \BU_{-w_\Delta(\alpha_1)}^\# \cdots
\BU_{-w_\Delta(\alpha_{i-1})}^\# \BU_{-w_\Delta(\alpha_{i+1})}^\#\cdots
 \BU_{-w_\Delta(\alpha_r)}^\#\}.$$
Let $\bar{X}'_i(w)=X'(w)\cup X'_i(w)\subset \BU_0$ (a disjoint union).
Let $X(w\cup w_i)= X(w)\cup X(w_i)$ (a disjoint union).
This is a partial compactification of $X(w)$, with divisor a union of
disconnected irreducible components and $X'_i(w)$ is obtained by removing
some of these components, as shows the following Proposition.

\begin{prop}
\label{partialcompactification}
We have a commutative diagram
$$\xymatrix{
\bar{X}'_i(w)\ar@{^{(}->}[r]^-{u\mapsto uw_\Delta\BB_0} &
 X(w\cup w_i) \\
X'_i(w) \ar@{^{(}->}[u] &
 X(w_i)\ar@{^{(}->}[u] \\
V_I \times X'_{\BL_I}(w_i)
 \ar[u]_\sim^{(v,u)\mapsto vw_\Delta w_I u w_I w_\Delta}
 \ar[d]^\sim_{(v,u)\mapsto (vw_\Delta w_I V_I,uw_I\BB_I)}\\
P_Iw_\Delta P_I/V_I\times_{L_I} X_{\BL_I}(w_i) \ar@{^{(}->}[r] &
 G/V_I\times_{L_I}X_{\BL_I}(w_i)\ar[uu]^\sim_{(gV_I,h\BB_i)\mapsto gh\BB_0}
}$$
\end{prop}

\begin{proof}
The map
$G/V_I\times_{L_I}X_{\BL_I}(w_i)\to
X(w_i)$
is an isomorphism \cite[Lemma 3]{Lu1}.
The map
$X'_{\BL_I}(w_i)\to X_{\BL_I}(w_i),\ u\mapsto uw_I\BB_I$
is an isomorphism by Theorem \ref{Lusztig} (b).
We are left with proving that the map
$$V_I \times X'_{\BL_I}(w_i)\to X'_i(w),\
(v,u)\mapsto vw_\Delta w_I u w_I w_\Delta$$
is an isomorphism.
Let $u=u_1u_2$ with $u_1\in \BV_I$ and $u_2\in \BU_I$. We have
$u^{-1}F(u)=u_2^{-1}F(u_2)v$ where
$v=F(u_2)^{-1}u_1^{-1}F(u_1)F(u_2)\in \BV_I$. So, $u\in X'_i(w)$ if and only
if $u_1\in V_I$ and $u_2\in X'_i(w)\cap \BU_I=
w_Iw_\Delta X'_{\BL_I}(w)w_\Delta w_I$ and we are done.
\end{proof}

\subsubsection{}
Let us describe now the variety $X''(w)$. Given $q'$ a power of $p$, we denote
by $\CL_{q'} : \BA^1 \to \BA^1,\ x \mapsto x^{q'} - x$ the Lang map.
This is an \'etale
Galois covering with group $\BF_{q'}$ ({\it Artin-Schreier covering}).
Given $\alpha \in [\Delta/\phi]$, we set $q_\alpha^*=1$ and given
$1 \le j \le d_\alpha-1$, we define inductively
$q_{\phi^j(\alpha)}^* = q_{\phi^{j-1}(\alpha)}^\circ q_{\phi^{j-1}(\alpha)}^*$.
We define
\begin{align*}
\gamma : \prod_{i=1}^r \CL_{q_{\alpha_i}}^{-1}(\BG_m) &\to X''(w) \\
(\xi_1,\dots,\xi_r) &\mapsto \prod_{i=1}^r \bigl(\prod_{j=0}^{d_{\alpha_i}-1}
x_{\phi^j(\alpha_i)} (\xi_i^{q_{\phi^j(\alpha_i)}^*}) \bigr) D(\BU_0)
\end{align*}
Note that the group $\prod_{i=1}^r \BF_{q_{\alpha_i}} \simeq U_0^F/D(\BU_0)^F$
acts on $\prod_{i=1}^r \CL_{q_{\alpha_i}}^{-1}(\BG_m)$ by addition on each
component.
Finally, define
$\CL : \prod_{i=1}^r \CL_{q_{\alpha_i}}^{-1}(\BG_m) \to (\BG_m)^r$,
$(\xi_1,\dots,\xi_r) \mapsto
(\CL_{q_{\alpha_1}}(\xi_1),\dots, \CL_{q_{\alpha_r}}(\xi_r))$.
This is the quotient map by $\prod_{i=1}^r \BF_{q_{\alpha_i}}$.

\smallskip
The next proposition is immediately checked:

\begin{prop}\label{prop isomorphisme}
The map $\gamma$ is an isomorphism of varieties. Through the isomorphism
$\prod_{i=1}^r \BF_{q_{\alpha_i}} \iso U_0/D(\BU_0)^F$, it is
$U_0/D(\BU_0)^F$-equivariant. Moreover, we have a commutative diagram
$$\xymatrix{
\prod_{i=1}^r \CL_{q_{\alpha_i}}^{-1}(\BG_m) \ar[d]_{\CL}
 \ar[rr]^{\gamma}_\sim && X''(w) \ar[d]^{\can} \\
(\BG_m)^r \ar[rr]_{\sim}^{\can} && \BU_0^F \backslash X'(w)
}$$
\end{prop}

\subsection{The variety $D(\BU_0)^F\backslash Y(\dot{w})$}
\label{secquo}
We describe the abelian covering $D(\BU_0)^F\backslash Y(\dot{w})$ of the
variety
$U_0\backslash X(w)\simeq \BG_m^r$ by the group $(U_0/D(\BU_0)^F)\times T$.

\subsubsection{}
Composing the isomorphism (b) of Theorem \ref{Lusztig}
and those of Propositions
\ref{quotient par DU} and \ref{prop isomorphisme}, we get an isomorphism
$D(\BU_0)^F\backslash X(w)\iso \prod_{i=1}^r \CL_{q_{\alpha_i}}^{-1}(\BG_m)$
and a commutative diagram whose squares are cartesian
\begin{equation}
\label{cartesien}
\xymatrix{
Y(\dot{w}) \ar[rr] \ar[d]_{\pi_w} && D(\BU_0)^F\backslash Y(\dot{w}) \ar[rr]
 \ar[d]_{\pi_w'} &&
 \BU_0^F \backslash Y(\dot{w}) \ar[d]_{\pi_w''} \\
X(w) \ar[rr] && \prod_{i=1}^r \CL_{q_{\alpha_i}}^{-1}(\BG_m)
 \ar[rr]_{\CL} && (\BG_m)^r
}
\end{equation}
Here, $\pi_w$, $\pi_w'$ and $\pi_w''$ are the quotient maps by the right
action of $\BT_0^{wF}$.

\subsubsection{}
\label{secGaloiscover}
Let $i : \hat{\BG} \to \BG$ be the semisimple simply-connected 
covering of the derived group of $\BG$. There exists 
a unique isogeny on $\hat{\BG}$ (still denoted by $F$) such that 
$i \circ F = F \circ i$. Let $Y^\circ(\dot{w})$ denote the image 
of the Deligne-Lusztig variety $Y_{\hat{\BG}}(\dot{w})$ through $i$. 
By Proposition \ref{Yestconnexe} and Theorem \ref{Lusztig}(c),
$Y^\circ(\dot{w})$ is connected. Its stabilizer in $G$ contains $i(\hat{\BG}^F)$, 
so in particular it is stabilized by $U_0$. 
It is also $F^d$-stable. Let $H$ be the stabilizer of $Y^\circ(\dot{w})$ 
in $T$ (we have $H=i(i^{-1}(\BT_0)^{wF})$).

We have a canonical $G\times (T\rtimes \langle F^d\rangle)^\oppose$-equivariant
isomorphism $Y^\circ(\dot{w})\times_H T\iso Y(\dot{w})$.

\smallskip
Recall that the tame fundamental group of $(\BG_m)^r$ is the $r$-th power of
the tame fundamental group of $\BG_m$.
There exists positive integers $m_1,\ldots,m_r$ dividing
$|H|$ and an \'etale Galois covering
$\pi'' : (\BG_m)^r \to U_0 \backslash Y^\circ(\dot{w})$ with Galois group
$N$ satisfying the following properties:
\begin{itemize}
\item $\pi_w'' \circ \pi'' : (\BG_m)^r \to (\BG_m)^r$ sends $(t_1,\dots,t_r)$
to $(t_1^{m_1},\dots,t_r^{m_r})$
\item
The restriction $\phi_i:\mu_{m_i}\to H$ of
the canonical map $\prod_{i=1}^r\mu_{m_i}\to H$ to the $i$-th factor is
injective.
\end{itemize}
So, $N$ is a subgroup of $\prod_{i=1}^r\mu_{m_i}$ and we have a canonical
isomorphism $(\prod_{i=1}^r\mu_{m_i})/ N \iso H$.
Moreover, $\pi''$ induces an isomorphism
$(\BG_m)^r / N \iso U_0 \backslash Y^\circ(\dot{w})$.

\medskip
Let us recall some constructions related to tori and their characters.
Let $Y(\BT_0)$ be the cocharacter group of $\BT_0$.
Let $c$ be a positive integer divisible by $d$
such that $(wF)^c=F^c$.
Let $\zeta$ be a generator of $\BF_{q^c}^\times$.
We consider the surjective morphism of groups
$$N_w:Y(\BT_0)\to T,\ \lambda\mapsto N_{F^c/wF}(\lambda)(\zeta).$$
We put $\beta_i^\vee=s_{\alpha_1}\cdots s_{\alpha_{i-1}}(\alpha_i^\vee)$.

\begin{prop}
\label{calculmi}
The subgroup $\phi_i(\mu_{m_i})$ of $T$ is generated by
$N_w(\beta_i^\vee)$.
\end{prop}

\begin{proof}
We identify $X(w)$ with $X'(w)$ and
$U_0\backslash X(w)$ with $\BG_m^r$ via the canonical isomorphisms.

Let $T_i=\BT^{w_iF}$.
Let $j_i:X(w)\to
X(w\cup w_i)$ and $k_i:X(w_i)\to X(w\cup w_i)$
be the canonical immersions.
Let $\theta\in\Irr(KT)$.
Let $\CF_\theta=(\pi_{w*}K)\otimes_{KT}K_\theta$.
Let $\CG_\theta=k_i^*j_{i*}\CF_\theta$.
From \cite[Proposition 7.3]{BonRou}, we deduce the following:
\begin{itemize}
\item[(case 1)]
if $\theta(N_w(\beta_i^\vee))\not=1$,
then, $\CG_\theta=0$
\item[(case 2)]
otherwise, there is a character
$\theta_i$ of $T_i$ such that $\CG_\theta\simeq \CF_{\theta_i}$, where
$\CF_{\theta_i}=(\pi_{w_i*}K)\otimes_{KT_i}K_{\theta_i}$.
\end{itemize}

Assume we are in case 2. Let $q_i:X(w_i)\to U_0\backslash X(w_i)$ be the
quotient map. Then,
$(q_{i*}\CF_{\theta_i})^{U_0}$ is a non-zero
locally constant sheaf $\CL_i$ of constant rank
on $U_0\backslash X(w_i)$. Via the
isomorphisms of Proposition \ref{partialcompactification}, its restriction
to $U_0\backslash X'_i(w)\simeq \BG_m^{i-1}\times \BG_m^{r-i}$ is a non-zero
locally constant sheaf. On the other hand, in case 1,
then the restriction of $(q_{i*}\CG_\theta)^{U_0}$ to
$\BG_m^{i-1}\times \BG_m^{r-i}$ is $0$.

Let $\bar{\CF}_\theta$ be the rank $1$ locally constant sheaf on
$\BG_m^r$ corresponding, via the covering $\pi_w''$, to the character
$\theta$. Let $q:X'(w)\to\BG_m^r$ be the canonical map.
Then, $\bar{\CF}_\theta\simeq (q_*\CF_\theta)^{U_0}$.
Denote by $j'_i:\BG_m^r\to \BG_m^{i-1}\times\BA^1\times \BG_m^{r-i}$ and
$k'_i:\BG_m^{i-1}\times \BG_m^{r-i}\to
\BG_m^{i-1}\times\BA^1\times \BG_m^{r-i}$ the canonical immersions.
Via base change, we deduce that
$k_i^{\prime *}j'_{i*}\bar{\CF}_\theta$ is $0$ in case 1 and
non-zero in case 2.
This means that the restriction of $\theta$ to $\mu_{m_i}$ is trivial in
case 1 but not in case 2.

We summarize the constructions in the following commutative diagram
$$\xymatrix{
Y(\dot{w})\ar[r]^-{\pi_w}\ar[d] & X(w)\simeq X'(w)\ar@{^{(}->}[r]\ar[d]_{q} &
 \bar{X}'_i(w)\ar[d] & X'_i(w)\ar@{_{(}->}[l]\ar[d] \\
U_0\backslash Y(\dot{w})\ar[r] & U_0\backslash X(w)\simeq
 \BG_m^r\ar@{^{(}->}[r]_-{j'_i} &
 \BG_m^{i-1}\times\BA^1\times \BG_m^{r-i} &
 \BG_m^{i-1}\times \BG_m^{r-i} \ar@{_{(}->}[l]^-{k'_i}
}$$
\end{proof}

Let $Y'=(D(\BU_0)^F \backslash Y^\circ(\dot{w}))
\times_{\BU_0^F \backslash Y^\circ(\dot{w})} (\BG_m)^r$. The cartesian square
(\ref{cartesien}) gives an isomorphism
$Y' \iso \prod_{i=1}^r Y_{q_{\alpha_i},m_i}$
where, given $q'$ is a power of $p$ and $s$ a positive integer, we set
$$Y_{q',s} = \{(\xi,t) \in \BA^1 \times \BG_m~|~ \xi^{q'} - \xi = t^s\}.$$
This variety has an action of $\BF_{q'}$ by addition on the first coordinate
and an action of $\mu_s$ by multiplication on the second coordinate.
We consider the quotient maps
$\tau_{q',s} : Y_{q',s} \to \BG_m$, $(\xi,t) \mapsto t$ and
$\rho_{q',s} : Y_{q',s} \to \BA^1-\BA^1(\BF_{q'})$, $(\xi,t) \mapsto \xi$.

Let $\tau = \prod_{i=1}^r \tau_{q_{\alpha_i},m_i} :
\prod_{i=1}^r Y_{q_{\alpha_i},m_i} \to (\BG_m)^r$,
$\pi' : \prod_{i=1}^r Y_{q_{\alpha_i},m_i} \to D(\BU_0)^F \backslash
 Y^\circ(\dot{w})$ the canonical map,
$\rho' = \prod_{i=1}^r \rho_{q_{\alpha_i},m_i} :
 \prod_{i=1}^r Y_{q_{\alpha_i},m_i}
 \to \prod_{i=1}^r \CL_{q_{\alpha_i}}^{-1}(\BG_m)$ and
$\rho''=\pi_w'' \circ \pi''$.
We have a commutative diagram all of whose squares are cartesian
%\equat\label{diagramme important}
$$\xymatrix{
\prod_{i=1}^r Y_{q_{\alpha_i},m_i} \ar[rr]^{\tau}
 \ar[d]_{\pi'} \ar@/_2cm/[dd]_{\rho'} &&
(\BG_m)^r \ar[d]_{\pi''} \ar@/^1.5cm/[dd]^{\rho''}\\
D(\BU_0)^F \backslash Y^\circ(\dot{w}) \ar[rr] \ar[d]_{\pi_w'} &&
\BU_0^F \backslash Y^\circ(\dot{w}) \ar[d]_{\pi_w''} \\
\prod_{i=1}^r \CL_{q_{\alpha_i}}^{-1}(\BG_m)
 \ar[rr]^{\CL} && (\BG_m)^r
}$$%\endequat

%\subsubsection{}
%\label{RGammatensor}
%We deduce from \S \ref{secGaloiscover} isomorphisms in $D^b(\CO T)$:
%$$R\Gamma_c(D(\BU_0)^F\backslash Y(\dot{w}),\CO)\iso
%R\Gamma_c(\prod_{i=1}^r Y_{q_{\alpha_i},m_i},\CO)
%\otimes_{\CO(\prod_{i=1}^r\mu_{m_i})}^\BL \CO T$$
%$$
%\iso
%\left(R\Gamma_c(Y_{q_{\alpha_1},m_1},\CO)\otimes_{\CO\mu_{m_1}}\CO T\right)
%\otimes_{\CO T}^\BL\cdots \otimes_{\CO T}^\BL
% \left(R\Gamma_c(Y_{q_{\alpha_r},m_r},\CO)\otimes_{\CO\mu_{m_1}}
% \CO T\right).
%$$

\subsection{Cohomology of $D(\BU_0)^F\setminus Y(\dot{w})$}
In this subsection, we shall describe the action of $\CO T$ on
$R\Gamma_c(D(\BU_0)^F\setminus Y(\dot{w}),\CO)$
and the action of $F^d$ on its cohomology. For this, we
determine first the cohomology of the curves $Y_{q',s}$ in order 
to use the previous diagram (\S \ref{seccohYqd}). We also need 
a result on Koszul complexes (\S \ref{section koszul}).

\subsubsection{Cohomology of the curves $Y_{q',s}$}
\label{seccohYqd}

Let $s\in\BZ_{>0}$ be prime to $p$ and let $q'$ be a power of $p$.
In \S \ref{secGaloiscover}, we
introduced a smooth affine connected closed curve $Y_{q',s}$
in $\BA^1\times\BG_m$.
This is a variety defined over $\BF_{q'}$ and we denote by $F'$ the
corresponding Frobenius endomorphism.
Note that $Y_{q',s}$ is an open subvariety of
a variety considered by Laumon \cite[\S 3.2]{Lau}. The
variety $Y_{q',2}$ was considered by Lusztig \cite[p.18]{Lu3}.
We will describe the complex of cohomology for the finite group actions
and the cohomology for the additional action of the Frobenius endomorphism.

\smallskip
Let $g$ be a generator of $\mu_s$. We put
$Z=0\to \CO\mu_s\xrightarrow{g-1}\CO\mu_s \to 0$, a complex
of $\CO\mu_s$-modules, with non-zero terms in degrees $0$ and $1$. Note
that, up to isomorphism of complexes of $\CO\mu_s$-modules, $Z$ does not
depend on the choice of $g$.

\begin{lemma}
\label{cohYqd}
Let $\psi\in\Irr(\BF_{q'})$.
We have
$$e_\psi R\Gamma_c(Y_{q',s},\CO)\simeq
\begin{cases}
\CO \mu_s[-1] & \textrm{ if }\psi\not=1\\
Z[-1] & \textrm{ if }\psi=1
\end{cases}
$$
in $D^b(\CO \mu_s)$.
\end{lemma}

\begin{proof}
We put $Y=Y_{q',s}$ in the proof.
Since $Y$ is a smooth affine curve, it follows that
$H^i_c(Y,\CO)=0$ for $i\not=1,2$ and
$H^*_c(Y,\CO)$ is free over $\CO$. Note further that
the action of $\mu_s$ is free on $Y$, so
$R\Gamma_c(Y,\CO)$ is a perfect complex of $\CO\mu_s$-modules.
Choose $R\Gamma_c(Y,\CO)\in D^b(\CO(\BF_{q'}\times\mu_s))$ the
unique (up to isomorphism) bounded complex of projective modules
with no non-zero direct summand homotopy equivalent to $0$. Then,
$R\Gamma_c(Y,\CO)$ is 
a complex of finitely generated projective $\CO(\BF_{q'}\times\mu_s)$-modules
with non-zero terms only in degrees $1$ and $2$.
Since $Y$ is irreducible, we have $H^2_c(Y,\CO)\simeq
\CO$, with a trivial action of $\BF_{q'}\times\mu_s$.

Given $\psi\in\Irr(\BF_{q'})$, let
$C_\psi=e_\psi R\Gamma_c(Y,\CO)$.
We have
$R\Gamma_c(Y,\CO)=\bigoplus_{\psi\in\Irr(\BF_{q'})}
C_\psi$.
Assume $\psi\not=1$. The complex $C_\psi$ has non-zero homology only
in degree $1$ and the $\CO \mu_s$-module $M=H^1(C_\psi)$ is projective.
As in \cite[\S 3.4]{DeLu}, one sees that the character of $M$
vanishes outside $1$, hence it is a free $\CO\mu_s$-module.
We have $M^{\mu_s}\simeq e_\psi H^1_c(Y/\mu_s,\CO)\simeq\CO$.
It follows that $M$ is a free $\CO\mu_s$-module of rank $1$.

We have $C_1\simeq R\Gamma_c(Y/\BF_{q'},\CO)$.
So, $H^1(C_1)\simeq H^2(C_1)\simeq \CO$ with trivial action of
$\mu_s$. It follows that $C_1^2$ is a projective cover of $\CO$ and $C_1$
is isomorphic to
$0\to \CO\mu_s\xrightarrow{g-1}\CO\mu_s \to 0$.
\end{proof}

Let $\theta\in\Irr(\mu_s)$. Let $c_\theta$ be the smallest positive
integer such that $F^{\prime c_\theta}$ fixes $\theta$. We denote by
$\tilde{\theta}$ the extension of $\theta$ to
$\mu_s\rtimes\langle F^{\prime c_\theta}\rangle$ on which
$F^{\prime c_\theta}$ acts trivially and we put
$L_\theta=\Ind_{\mu_s\rtimes\langle F^{\prime c_\theta}\rangle}^{
\mu_s\rtimes\langle F'\rangle}\tilde{\theta}$.

Given $\lambda\in K^\times$, we denote by $K(\lambda)$ the 
vector space $K$ with action of $F'$ given by multiplication by $\lambda$.
Let $Q=\{\alpha\in\bar{\BF}_{p}^\times|\alpha^s\in\BF_{q'}^\times\}$.

\begin{lemma}
\label{FonYqd}
Fix $\psi\in\Irr(\BF_{q'})^\#$.
There is $\lambda:\Irr(\mu_s)^\#/F'\to \CO^\times$ such that
$$H^1_c(Y_{q',s},K)\simeq \bigl(K\BF_{q'}\otimes K(1)\bigr)\oplus
\bigoplus_{\substack{\alpha\in Q/\mu_s\\
\theta\in[\Irr(\mu_s)^\#/F']}}
\left(\alpha^s(\psi)\otimes L_\theta\otimes
K\left(\lambda(\theta)\cdot\sqrt[c_\theta]{
\theta(\alpha^{q^{\prime c_\theta}-1})}\right) \right)$$
as $K(\BF_{q'}\times\mu_s\rtimes\langle F'\rangle)$-modules.

Furthermore, 
we have $H^2_c(Y_{q',s},K)\simeq K(q')$.
\end{lemma}

\begin{proof}
Note first that the statement about $H^2_c(Y_{q',s},K)$ follows immediately
from Lemma \ref{cohYqd}.

We have an action of $Q$ on $Y_{q',s}$ extending the action of $\mu_s$:
an element $\alpha\in Q$ acts by
$(\xi,t)\mapsto (\alpha^s\xi,\alpha t)$.
This provides $Y_{q',s}$ with an action of
$(\BF_{q'}\rtimes Q)\rtimes\langle F'\rangle$. The action of
$Q$ on $\BF_{q'}$ is given by $\alpha:(\BF_{q'}\ni a\mapsto \alpha^s a)$.

Let $L=H^1_c(Y_{q',s},K)/H^1_c(Y_{q',s},K)^{\mu_s}$.
 Given $\theta\in\Irr(\mu_s)$, we put
$V_\theta=\Ind_{\BF_{q'}\times\mu_s}^{\BF_{q'}\rtimes Q}(\psi\otimes\theta)$
(this is independent of the choice of $\psi$ up to isomorphism).
Then, $L\simeq \bigoplus_{\theta\in\Irr(\mu_s)^\#}V_\theta$ is the
decomposition into irreducible $K(\BF_{q'}\rtimes Q)$-modules: it is
uniquely determined by its restriction to $\BF_{q'}\times\mu_s$ which
is described in Lemma \ref{cohYqd}.

Let $\theta\in\Irr(\mu_s)^\#$. We put
$L'_\theta=\Ind_{\BF_{q'}\times\mu_s\rtimes
 \langle F^{\prime c_\theta}\rangle}^{
\BF_{q'}\rtimes Q\rtimes\langle F'\rangle}(\psi\otimes\tilde{\theta})$.
This is an irreducible representation which extends
$\bigoplus_{i=0}^{c_\theta-1}V_{F^{\prime i}(\theta)}$.
It follows that there are scalars
$\lambda(\theta)\in \CO^\times$ such that
$$L\simeq \bigoplus_{\theta\in\Irr(\mu_s)^\#/F'}
L'_\theta\otimes K(\lambda(\theta)).$$
We have
$$\Res_{\BF_{q'}\times \mu_s\rtimes\langle F'\rangle}^{
\BF_{q'}\rtimes Q\rtimes\langle F'\rangle}L'_\theta\simeq
\Ind_{\BF_{q'}\times \mu_s\rtimes\langle F^{\prime c_\theta}\rangle}^{
\BF_{q'}\times \mu_s\rtimes\langle F'\rangle}
\Res_{\BF_{q'}\times \mu_s\rtimes\langle F^{\prime c_\theta}\rangle}^{
\BF_{q'}\rtimes Q\rtimes\langle F^{\prime c_\theta}\rangle}L''_\theta$$
where $L''_\theta=\Ind_{\BF_{q'}\times\mu_s\rtimes
 \langle F^{\prime c_\theta}\rangle}^{
\BF_{q'}\rtimes Q\rtimes\langle F^{\prime c_\theta}\rangle}
(\psi\otimes\tilde{\theta})$.
Now,
$$\Res_{\BF_{q'}\times \mu_s\rtimes\langle F^{\prime c_\theta}\rangle}^{
\BF_{q'}\rtimes Q\rtimes\langle F^{\prime c_\theta}\rangle}L''_\theta\simeq
\bigoplus_{\alpha\in Q/\mu_s} \alpha^s(\psi)\otimes\tilde{\theta}\otimes
K'(\theta(\alpha^{q^{\prime c_\theta}-1}))$$
where $K'(\lambda)$ is the one dimensional representation of
$\langle F^{\prime c_\theta}\rangle$ where $F^{\prime c_\theta}$ acts
by $\lambda$. The Lemma follows.
\end{proof}

\subsubsection{Koszul complexes}
\label{section koszul}
Let $H$ be a finite abelian group. Given 
$H'\le H$ a cyclic subgroup and $g$ a generator of $H'$, we put
$Z_H(H')=0\to \CO H \xrightarrow{g-1}\CO H\to 0$,
where the non-zero terms are in degrees $0$ and $1$. This is a complex
of $\CO H$-modules whose isomorphism class is independent of $g$.
Given $H_1,\ldots,H_m$ a collection of cyclic subgroups of $H$, we put
$Z_H(H_1,\ldots,H_m)=Z_H(H_1)\otimes_{\CO H}\cdots\otimes_{\CO H}Z_H(H_m)$,
a Koszul complex.

\begin{lemma}
\label{tensorZ}
Let $H_1,\ldots,H_m$ be finite cyclic subgroups of $H$.
Then, the cohomology of the complex 
$Z_H(H_1,\ldots,H_m)$ is free over $\CO$.
\end{lemma}

\begin{proof}
We prove the Lemma by induction on $m$.
Let $\bar{H}=H/H_n$ and $\bar{H}_i=H_i/(H_i\cap H_n)$.
We have a distinguished triangle in $D^b(\CO H)$
$$H^0 Z_H(H_n)\to Z_H(H_n)\to (H^1 Z_H(H_n))[-1]\rightsquigarrow.$$
We have an isomorphism
$\CO \bar{H}\iso H^1 Z_H(H_n)$ induced by the identity and an isomorphism
$\CO \bar{H}\iso H^0 Z_H(H_n)$ induced by $1\mapsto \sum_{h\in H_n}h$.
 We have an
isomorphism of complexes of $\CO\bar{H}$-modules
$$\CO\bar{H}\otimes_{\CO H}Z_H(H_1,\ldots,H_{n-1})\iso
Z_{\bar{H}}(\bar{H}_1,\ldots,\bar{H}_{n-1})$$
and we deduce that there is a distinguished triangle in $D^b(\CO H)$
$$Z_{\bar{H}}\to Z_H\to Z_{\bar{H}}[-1]\rightsquigarrow$$
where 
$Z_H=Z_H(H_1,\ldots,H_n)$ and
$Z_{\bar{H}}=Z_{\bar{H}}(\bar{H}_1,\ldots,\bar{H}_{n-1})$.
This induces a long exact sequence 
$$0\to H^0 Z_{\bar{H}}\to H^0 Z_H\to 0\to H^1 Z_{\bar{H}}\to H^1 Z_H\to
H^0Z_{\bar{H}}\to H^2 Z_{\bar{H}}\to\cdots$$
Note that $KZ_H(H')\simeq K(H/H')\oplus K(H/H')[-1]$, hence
$H^i(KZ_H)\simeq K(H/H_1\cdots H_n)^{n \choose i}$ and 
similarly
$H^i(KZ_{\bar{H}})\simeq K(H/H_1\cdots H_n)^{n-1\choose i}$.
We deduce that the connection maps $H^iZ_{\bar{H}}\to H^{i+2}Z_{\bar{H}}$
vanish over $K$. By induction, these $\CO$-modules are free over $\CO$, hence
these maps vanish over $\CO$ as well and it follows that 
$H^* Z_H$ is free over $\CO$.
\end{proof}

\subsubsection{Action of $T$, Conjecture \ref{conj2plus}}\label{RGammatensor}
We deduce from \S \ref{secGaloiscover} isomorphisms in $D^b(\CO T)$:
\begin{multline*}
R\Gamma_c(D(\BU_0)^F\backslash Y(\dot{w}),\CO)\iso
R\Gamma_c(\prod_{i=1}^r Y_{q_{\alpha_i},m_i},\CO)
\otimes_{\CO(\prod_{i=1}^r\mu_{m_i})}^\BL \CO T\\
\iso
\left(R\Gamma_c(Y_{q_{\alpha_1},m_1},\CO)\otimes_{\CO\mu_{m_1}}\CO T\right)
\otimes_{\CO T}^\BL\cdots \otimes_{\CO T}^\BL
 \left(R\Gamma_c(Y_{q_{\alpha_r},m_r},\CO)\otimes_{\CO\mu_{m_1}}
 \CO T\right).
\end{multline*}

%We put
%$\CW_j=\{q_{\alpha_{i_1}}\cdots q_{\alpha_{i_j}}\}_{
%1\le i_1<\cdots<i_j\le r}$. In the split case, we have $\CW_j=\{q^j\}$.

Let $a\in\CO T$. We put
$Z(a)=0\to \CO T\xrightarrow{a}\CO T\to 0$, with non-zero terms in
degrees $0$ and $1$.
Given $a_1,\ldots,a_m\in \CO T$, we define the Koszul complex
$Z(\{a_1,\ldots,a_m\})=Z(a_1)\otimes_{\CO T}\cdots\otimes_{\CO T}Z(a_m)$.
We also put $Z(\emptyset)=\CO T$. Given $g\in T$, we have
$Z(g-1)=Z_T(\langle g\rangle)$ with the notations of \S \ref{section koszul}.

 If $\psi$ is a linear character 
of $U_0$ trivial on $D(\BU_0)^F$, we put 
$I(\psi)=\{i \in \{1,2,\dots,r\}~|~\psi_{\alpha_i} \neq 1\}$. 
It is easily checked that, if $t \in T_0$, then $I(^t\psi)=I(\psi)$. 

\smallskip
%We fix a $G$-regular character $\psi$ of $U_0$. Given $I\subset\{1,\ldots,r\}$,
%we denote by $\psi_I'$ the character of $U_0$ whose
%restriction to $\BF_{q_{\alpha_i}}$ is equal to the restriction of $\psi$
%if $i\in I$ and is trivial otherwise.
%
%\smallskip

From the above isomorphism, Proposition \ref{calculmi}, and
Lemmas \ref{cohYqd} and \ref{tensorZ}, we deduce
\begin{lemma}
\label{cohYsurD}
Let $\psi$ be a linear character of $U_0$ trivial on $D(\BU_0)^F$. 
We have an isomorphism in $D^b(\CO T)$
$$e_\psi C[r]\simeq Z(\{N_w(\beta_i^\vee)-1\}_{i\not\in I(\psi)}).$$
The cohomology of $e_\psi C$ is free over $\CO$.
\end{lemma}

Conjecture \ref{conj2plus} for twisted Coxeter elements of Levi subgroups 
now follows easily from this lemma:

\begin{thm}
\label{coxeter}
Let $I$ be a $\phi$-stable subset of $\Delta$. Let
$[I/\phi]=\{\beta_1,\dots,\beta_{r'}\}$ and
$w'=s_{\beta_1}\cdots s_{\beta_{r'}}$.
Then, Conjecture \ref{conj2plus} holds for
$(\BG,w')$.
\end{thm}

\begin{proof}
Thanks to Remark \ref{reduceLevi}, we can and will assume that
$I=\Delta$. In other words, we may assume that $w'=w$. 
We have $Y_x(\dot{w})=\emptyset$ for $x\not=w_\Delta$ by
Theorem \ref{Lusztig} (a), so it is enough to prove Conjecture \ref{conj2}.

Let $\psi$ be a regular linear character of $U_0$. By Lemma \ref{cohYsurD},
we have
$e_\psi R\Gamma_c(Y(\dot{w}),\CO)\simeq Z(\emptyset)[-r] = \CO T[-r]$.
\end{proof}

\subsubsection{Action of $F^d$}\label{sec frobenius}
Consider the endomorphism of $\prod_{i=1}^r Y_{q_{\alpha_i},m_i}$ given
by elevation to the power $q_{\alpha_i}^{d/d_{\alpha_i}}$ on the $i$-th
component.  Then, the morphism
$\pi'$ is equivariant with respect to the action of that endomorphism and of
$F^d$ on $D(\BU_0)^F \backslash Y^\circ(\dot{w})$.

Let $\theta\in\Irr(T)$. Given $i\in\{1,\ldots,r\}$, we denote by $\theta_i$ 
the restriction of $\theta$ to $\phi_i(\mu_{m_i})=<N_w(\beta_i^\vee)>$. 
Let $I_\theta=\{i \in \{1,2,\dots,r\}~|~\theta_i \neq 1\}$. It is 
easily checked that $I_{\theta \circ F^d}=I_\theta$. 
Let $c_\theta$ be the smallest positive
integer such that $F^{dc_\theta}$ fixes $\theta$. We denote by
$\tilde{\theta}$ the extension of $\theta$ to
$T\rtimes\langle F^{dc_\theta}\rangle$ on which
$F^{dc_\theta}$ acts trivially and we put
$V_\theta=\Ind_{T\rtimes\langle F^{dc_\theta}\rangle}^{
T\rtimes\langle F^d\rangle}\tilde{\theta}$.
Given $\lambda\in K^\times$,
we denote by
$K(\lambda)$ the vector space $K$ with an action of $F^d$
given by multiplication by $\lambda$.

\begin{lemma}
\label{FsurCI}
Let $\psi$ be a linear character of $U_0$ trivial on $D(\BU_0)^F$. 
Then there is a map $\nu:\Irr(T)/F^d\to \CO^\times$ such that for any
$j\in\BZ$,
we have an isomorphism of $K(T\rtimes\langle F^d\rangle)$-modules
$$H^{r+j}(e_\psi KC)\simeq
\bigoplus_{\substack{\theta\in[\Irr(T)/F^d]\\ I_\theta \subset I(\psi)}}
\left(V_\theta\otimes K(q^{dj}\nu(\theta))
\right)^{\oplus {r-|I(\psi)|\choose j}}$$
for $0\le j\le r-|I(\psi)|$. Furthermore,
$H^{r+j}(e_\psi KC)=0$ if $j<0$ or $j>r-|I(\psi)|$.
\end{lemma}

\begin{proof}
We use Lemma \ref{FonYqd}. There are scalars $\lambda(i,\theta)\in \CO^\times$
such that for $i\in I(\psi)$, we have
$$e_{\psi_{\alpha_i}}H^1_c(Y_{q_{\alpha_i},m_i},K)\otimes_{K\mu_{m_i}}K T\simeq
\bigl(\bigoplus_{\substack{\theta\in [\Irr(T)/F^d]\\ \theta_i\not=1}}
 V_\theta\otimes K(\lambda(i,\theta))\bigr)
\oplus
\bigl(\bigoplus_{\substack{\theta\in [\Irr(T)/F^d]\\ \theta_i=1}} V_\theta\bigr).$$
Note that 
$e_{\psi_{\alpha_i}}H^j_c(Y_{q_{\alpha_i},m_i},K)=0$ for $j\not=1$.
On the other hand, we have
$$e_1 H^1_c(Y_{q_{\alpha_i},m_i},K)\otimes_{K\mu_{m_i}}K T\simeq
\bigoplus_{\substack{\theta\in[\Irr(T)/F^d]\\ \theta_i=1}}
V_\theta$$
and
$$e_1 H^2_c(Y_{q_{\alpha_i},m_i},K)\otimes_{K\mu_{m_i}}K T\simeq
\bigoplus_{\substack{\theta\in[\Irr(T)/F^d]\\ \theta_i=1}}
 V_\theta\otimes K(q^d).$$
Note that 
$e_1H^j_c(Y_{q_{\alpha_i},m_i},K)=0$ for $j\not=1,2$.

Given $\theta\in \Irr(T)$, we have $V_\theta\otimes_{KT}V_\theta\simeq
V_\theta$.
So,
$$H^{r+j}(e_\psi KC)\simeq
\left(\bigoplus_{\substack{\theta\in[\Irr(T)/F^d]\\ I_\theta \subset I(\psi)}}
V_\theta\otimes K\left(q^{dj}\prod_{i\in I_\theta}
 \lambda(i,\theta)\right)\right)^{\oplus{r-|I(\psi)|\choose j}}$$
for $0\le j\le r-|I(\psi)|$ and
$H^{r+j}(e_\psi KC)=0$ otherwise. The result follows.
\end{proof}

From now on, and until the end of this section, we fix a regular linear character 
$\psi$ of $U_0$. We put 
$$e(\psi)=\sum_{\substack{I \subset \Delta \\ \phi(I)=I}} e_{\tilde{\psi}_I}.$$
For the definition of $\tilde{\psi}_I$, see \S \ref{psi I}. 

\begin{rem}\label{orbites}
If the center of $\BG$ is connected, then 
$\{\tilde{\psi}_I~|~I \subset \Delta,~\phi(I)=I\}$ is a set of 
representatives of $T_0$-orbits of linear characters of $U_0$ which are 
trivial on $D(\BU_0)^F$. This follows easily from \cite[Theorem 2.4]{DiLeMi1}
and from the fact that the centre of any Levi subgroup of $\BG$ is also 
connected \cite[Lemma 1.4]{DiLeMi1}.
\end{rem}

\begin{prop}
\label{FsurC}
There is a map $\nu:\Irr(T)/F^d\to \CO^\times$ such that for all $j$,
we have an isomorphism of $K(T\rtimes\langle F^d\rangle)$-modules
$$H^{r+j}(e(\psi)KC)\simeq
\mathop{\oplus}_{\substack{I \subset \{1,2,\dots,r\}\\ |I|\le r-j}} \bigl(
\bigoplus_{\substack{\theta\in[\Irr(T)/F^d]\\ I_\theta \subset I}}
\left(V_\theta\otimes K(q^{dj}\nu(\theta))
\right)^{\oplus {r-|I|\choose j}}\bigr).$$
\end{prop}

\begin{proof}
This follows easily from the proof of Lemma \ref{FsurCI}, the scalars 
$\lambda(i,\theta)$ (which are defined whenever $i \in I(\psi) \cap I_\theta$) 
depending only on $\theta$ and $\psi_{\alpha_i}$.
\end{proof}

\subsubsection{Endomorphism algebra}
Let $E=R\End^\bullet_{\CO T}(e(\psi)C)$, an object of $D^b(\CO T)$.

\begin{prop}
\label{descE}
We have an isomorphism in $D^b(\CO T)$
$$E\simeq \bigoplus_{\substack{I\subset\{1,\ldots,r\} \\ j \in \BZ}}
 Z(\{N_w(\beta_i^\vee)-1\}_{i\in I})[j]^{m_{I,j}}$$
for some non-negative integers $m_{I,j}$.
The cohomology of $E$ is free over $\CO$ and
$H^j(E)=0$ for $|j|>r$.

Let $j\in\BZ$.
We have an isomorphism of $K(T\rtimes\langle F^d\rangle)$-modules
$$H^j(KE)\simeq \bigoplus_{\theta\in[\Irr(T)/F^d]} \left(V_\theta\otimes
 K({q^{dj}})\right)^{m(\theta,j)}$$
for some integers $m(\theta,j)$. We have $m(\theta,j)=0$ if
$r-|I_\theta|<|j|$.
\end{prop}

\begin{proof}
Given $a\in\CO T$, we have isomorphisms
$$\Hom^\bullet_{\CO T}(Z(a),\CO T)\simeq Z(a)[1] \textrm{ and }
Z(a)\otimes_{\CO H}Z(a)\simeq Z(a)\oplus Z(a)[-1].$$
We deduce from Lemma \ref{cohYsurD} that
$e(\psi)C\simeq \bigoplus_{I\subset\{1,\ldots,r\}}
Z(\{N_w(\beta_i^\vee)-1\}_{i\in I})[-r]$. Hence,
$$E\simeq \bigoplus_{I,J\subset\{1,\ldots,r\}}\bigoplus_{k=0}^{|I\cap J|}
Z(\{N_w(\beta_i^\vee)-1\}_{i\in I\cup J})[|J|-k]^{{|I\cap J|\choose k}}.$$
The freeness of the cohomology follows from
Lemma \ref{tensorZ}.
The second assertion is a direct consequence of
Proposition \ref{FsurC}.
\end{proof}

From now on, and until the end of this section, we assume that $\BG$ is quasi-simple. 
Let $h$ be the Coxeter number of $\BG$ relative to $F$
(cf \cite[\S 1.13]{Lu2}). In other words, $h=|W^{wF}|$. 
Let $\nu$ be a non-negative integer with $\ell^\nu\equiv 1\pmod h$.
We put $\tilde{F}=F^{d\ell^\nu}$.

\begin{cor}
\label{specE}
Assume that $\BG$ is quasi-simple and that $\ell{\not|}h$. Then, the canonical map
$k\otimes_{\CO}H^0(E)^{\tilde{F}}\to H^0(k\otimes_{\CO}E)^{\tilde{F}}$ is an
isomorphism.

Assume furthermore that for all $\theta\in\Irr(T)$ and all
$j\in\{1,\ldots,r-|I_\theta|\}$, then $q^{djc_\theta}\not=1\pmod \ell$.
Then, $H^i(E)^{\tilde{F}}=H^i(k\otimes_{\CO}E)^{\tilde{F}}=0$ for $i\not=0$.
\end{cor}

\begin{proof}
Note that the canonical map $kH^i(E)\to H^i(kE)$ is an
isomorphism since $H^i(E)$ is free over $\CO$. Now,
$H^i(E)^{\tilde{F}}$ coincides with the generalized eigenspace of
$\tilde{F}$ for the
eigenvalue $1$. The eigenvalues of $\tilde{F}$ on $H^0(E)$ are 
$h$-th roots of unity, hence their reductions modulo $\ell$ remain
distinct. It follows that the generalized $1$-eigenspace of $\tilde{F}$ on
$kH^0(E)$ is the image of the generalized $1$-eigenspace on $H^0(E)$.

Also, the eigenvalues of $\tilde{F}$ on $H^j(E)$ and $H^j(kE)$
are of the form $q^{\ell^\nu dj}\zeta_{c_\theta}$, where
$\zeta_{c_{\theta}}^{c_{\theta}}=1$ and the result follows.
\end{proof}

%\subsection{Disjunction}
%\subsubsection{}
%Let $M$ be the sum of induced Gelfand-Graev representations as in
%\S \ref{secgeneration}. Let $A_M=\End_{\CO G}(M)$,
%$C_M=R\Hom^\bullet_{\CO G}(M,R\Gamma_c(Y(\dot{w}),\CO))$,
%and $E_M=R\End_{\CO \BT_0^{wF}}^\bullet(C_M)$.
%We have
%$$R\Hom^\bullet_{\CO G}(M,R\Gamma_c(Y(\dot{w}),\CO))\simeq
%\bigoplus_{I\subset\Delta}e_{\tilde{\psi}_I}
% R\Gamma_c(D(\BU_0)^F\backslash Y(\dot{w}),\CO)$$
%Since $R\Gamma_c(D(\BU_0)^F\backslash Y(\dot{w}),\CO)\simeq
%R\Gamma_c(D(\BU_0)^F\backslash Y^\circ(\dot{w}),\CO)\otimes_{\CO H}\CO \BT_0^{wF}$,
%we deduce from Proposition \ref{descE} and Corollary \ref{specE}:
%
%\begin{lemma}
%\label{endGG}
%The $\CO$-module $H^*(E_M)$ is free over $\CO$ and $H^i(E_M)=0$ for $|i|>r$.
%
%If ???, then
%$H^i(E_M)^{\tilde{F}}=0$ for $i\not=0$ and
%the canonical map
%$k\otimes_{\CO}H^0(E_M)^{\tilde{F}}\to H^0(k\otimes_{\CO}E_M)^{\tilde{F}}$ is
%an isomorphism.
%\end{lemma}
%

\section{Groups of type $A$}

\begin{quotation}
\noindent{\bf Hypothesis:} {\it In this section, and only in this section, 
we assume that $\BG$ is of type $A_{n-1}$ (for some non-zero natural number $n$) 
and that $F$ is a split Frobenius endomorphism of $\BG$ (i.e. $d=1$). 
We keep the notation of Section \ref{section coxeter}, i.e. $w$ is a 
Coxeter element of $W$.}
\end{quotation}

Note in particular that $r=n-1$, $W \simeq \GS_n$, $w$ is a 
cycle of length $n$, and $h=n$. 

\subsection{A progenerator for $\CO G$}
\label{secgeneration}
%Let 
%$$M=\mathop{\oplus}_{\psi \in [\Hom(U_0,\CO^\times)/T_0]} \Gamma_\psi.$$
%Since $\BG$ is of type $A$, we have $D(\BU_0)^F=D(U_0)$, so every 
%linear character of $U_0$ is trivial on $D(\BU_0)^F$. Therefore, 
%every $\Gamma_\psi$ in the above direct sum is a Harish-Chandra 
%induced of a Gelfand-Graev module (see \ref{hc induction}). 
%
%\begin{rem}
%Assume in this remark that the centre 
%of $\BG$ is connected. Let $\psi$ be a regular linear character 
%of $U_0$. Then $\{\tilde{\psi}_I~|~I \subset S\}$ is a set of 
%representatives of $T_0$-orbits in $\Hom(U_0,\CO^\times)$.
%\end{rem}
%
The following result is probably classical.

\begin{thm}
\label{progenerator}
$\Ind_{D(\BU_0)^F}^G \CO$ is a progenerator for $\CO G$.
\end{thm}

\begin{proof}
Consider the situation of \S \ref{section reductions}. 
Let $\tilde{\BU}_0=i(\BU_0)$. Then $i$ induces an isomorphism 
$D(\BU_0)^F\iso D(\tilde{\BU}_0)^F$. So,
$$\Ind_{D(\tilde{\BU}_0)^F}^{\tilde{G}} \CO \simeq
\Ind_G^{\tilde{G}}\bigl(\Ind_{D(\BU_0)^F}^G \CO\bigr)$$
and
$$\Res_G^{\tilde{G}} \Ind_{D(\tilde{\BU}_0)^F}^{\tilde{G}} \CO 
\simeq \Bigl(\Ind_{D(\BU_0)^F}^G \CO\Bigr)^{\oplus r},$$
where $r=|\tilde{G}/i(G)|$.
It follows that $\Ind_{D(\BU_0)^F}^G \CO$ is a progenerator for $\CO G$
if and only 
if $\Ind_{D(\tilde{\BU}_0)^F}^{\tilde{G}} \CO$ is a progenerator 
for $\CO\tilde{G}$. This implies that we only need to prove the theorem 
for $\BG=\GL_n$, which we assume for the rest of the proof.
Let $M_I=\Ind_{D(\BU_I)^F}^{L_I} \CO$. 

Let $S$ be a simple cuspidal $\CO G$-module. Then,
there is an $\CO G$-module $\tilde{S}$, free over $\CO$, with
$k\tilde{S}\simeq S$ and such that $K\tilde{S}$ is cuspidal 
(see \cite[4.15 and 5.23]{dipper}). 
Now, there exists a regular linear character $\psi$ of $U_0$ 
such that $\Hom_{KG}(K\Gamma_\psi,K\tilde{S})\not=0$, so
$\Hom_{\CO G}(\Gamma_\psi,\tilde{S})\not=0$, since
$\Gamma_\psi$ is projective. Hence $\Hom_{\CO G}(\Gamma_\psi,S)\not=0$. 
In particular, $\Hom_{\CO G}(M_\Delta,S)\not=0$, because $\Gamma_\psi$ 
is a direct summand of $M_\Delta$. 

Let now $S$ be an arbitrary simple $\CO G$-module. Then, there is
$I\subset\Delta$ and a simple cuspidal $\CO L_I$-module $S'$ such that
$\Hom_{\CO G}(\CR_{\BL_I\subset\BP_I}^\BG S',S)\not=0$ and we deduce
from the study above and the exactness of the Harish-Chandra induction 
that $\Hom_{\CO G}(\CR_{\BL_I\subset\BP_I}^\BG M_I,S)\not=0$.
It follows again that $\Hom_{\CO G}(M_\Delta,S)\not=0$ for any simple $\CO G$-module $S$ 
because $\CR_{\BL_I\subset\BP_I}^\BG M_I$ is a direct summand of $M_\Delta$. 
So we are done.
\end{proof}

\begin{cor}
\label{A centre connexe}
Let $\psi$ be a regular linear character of $U_0$. 
If the centre of $\BG$ is connected, then $\CO G e(\psi)$ is 
a progenerator for $\CO G$.
\end{cor}

\begin{proof}
This follows from Theorem \ref{progenerator} and Remark \ref{orbites}.
\end{proof}

%\begin{rem}
%It is of course enough, in the definition of $M$, to sum over
%the characters $\psi$ modulo $T_0$-conjugation. In particular, if
%the center of $\BG$ is connected, then one $\psi$ suffices.
%\end{rem}

\subsection{Description over $K$} 
We identify $T$ with $T_w$ as in \S \ref{modelsRTG}.
From now on, and until the end of this paper, we assume that $\ell$ divides $|G|$ 
but does not divide $|G/T|$. In other words, we assume that the order of $q$ 
modulo $\ell$ is equal to $n$. Note in particular that $\ell > n$. 
Let $S$ be a Sylow $\ell$-subgroup of $T$ (note that $S$ is a Sylow $\ell$-subgroup 
of $G$). Let $b$ be the sum of all 
block idempotents of $\CO G$ having $S$ as a defect group.
Since $C_\BG(s)=\BT_w$ for every non-trivial 
$\ell$-element $s$ of $T$, $b$ is the sum of all block idempotents of 
$\CO G$ with non-zero defect.

\begin{lemma}
\label{EndK}
The canonical map $bKG\to \End_{KT}^\bullet(KC)^{\tilde{F}}$ is an isomorphism
in $D^b(KG\otimes (KG)^\oppose)$.
\end{lemma}

\begin{proof}
Let $\zeta,\zeta'\in\Irr_K(T)$ be two characters which are not conjugate
under the action of $N_G(T)$. Then,
$\Hom_{KG}(H^*_c(Y(\dot{w}),K)\otimes_{KT}\zeta,
H_c^*(Y(\dot{w}),K)\otimes_{KT}\zeta')=0$ \cite[Proposition 13.3]{DiMi2}.

Let $\zeta\in\Irr(T)$.
Let $s$ be a semi-simple element
of the group dual to $\BG$ whose class is dual to that of $\zeta$.
Let $\BL$ be an $F$-stable Levi subgroup of
$\BG$ corresponding to the centralizer of $s$ and containing $\BT_w$.
There is a decomposition 
$n=n_1n_2$ in positive integers such that $L=\GL_{n_1}(\BF_{q^{n_2}})$ and
$\BT_w$ is a Coxeter torus of $\BL$. We denote by $Y_L$ the Deligne-Lusztig
variety associated to a Coxeter element for $\BL$. Let $\BP$ be a parabolic
subgroup of $\BG$ with Levi complement $\BL$.

The character $\zeta$ is the restriction of a character
$\hat{\zeta}\in\Irr_K(L)$
and we have an isomorphism of $KL$-modules
$$H^j_c(Y_L,K)\otimes_{KT}\zeta\simeq \hat{\zeta}\otimes_K
H^j_c(Y_L/T,K).$$
By \cite[Theorem 6.1]{Lu2}, we have
$$\Hom_{KL}(H^i_c(Y_L/T,K),H_c^{i'}(Y_L/T,K))=0
\textrm{ for } i\not=i',$$
hence
$$\Hom_{KL}(H^i_c(Y_L,K)\otimes_{KT}\zeta,H_c^{i'}(Y_L,K)\otimes_{KT}\zeta)=0
\textrm{ for } i\not=i'.$$
Furthermore, if not zero, then $H^j_c(Y_L/T,K)$ is an irreducible representation
of $KL$ \cite{Lu2}.
By \cite[Theorem B]{BonRou}, the $KG$-modules
$\CR_{\BL\subset \BP}^\BG(\hat{\zeta}\otimes_K H^i_c(Y_L/T,K))$ are
irreducible and distinct, when $i$ varies. Since
$$\CR_{\BL\subset \BP}^\BG(R\Gamma_c(Y_L,K))\simeq
R\Gamma_c(Y(\dot{w}),K),$$
it follows that the $H^i_c(Y(\dot{w}),K)\otimes_K \zeta$ are
irreducible and distinct, when $i$ varies.

Note that $H_c^*(Y(\dot{w}),K)\otimes_{KT}\zeta$ depends
on $\zeta$ only up to $N_G(T)$-conjugacy.
Since the action of $\tilde{F}$ on $T$ coincides with that of $w$, it follows
that, if $V$ is an irreducible $KG$-module, then $V \otimes_{KG} H^*_c(Y(\dot{w}),K)$ 
is irreducible under the right action of
$T\rtimes \langle\tilde{F}\rangle$.
The result follows now from Lemma \ref{local}.
\end{proof}

\begin{rem}
One needs actually a weaker statement than the disjunction of the
cohomology groups for the $G$-action. We only need to know that the
part of the Lefschetz character of the cohomology of $Y(\dot{w})$ corresponding
to a fixed eigenvalue of $F$ is an irreducible representation of $G^F$.
This can be deduced, in the unipotent case, from \cite{DiMi1}.
\end{rem}

\subsection{Determination of the $(T\rtimes F)$-endomorphisms} 

\begin{prop}
\label{EndLambda}
For $R=\CO$ or $k$, the
canonical map
$bR G\to R\End_{R T}^\bullet(RC)^{\tilde{F}}$ 
is an isomorphism in $D^b(RG\otimes (RG)^\oppose)$.
\end{prop}

\begin{proof}
By Proposition \ref{isobetweenY}, it is enough to prove the Proposition for the
simply connected covering of $[\BG,\BG]$, which is isomorphic to $\SL_n$.
Using the embedding
$\SL_n\to\GL_n$, we have a further reduction to the case $\BG=\GL_n$.

We choose for $\BT_0$ the diagonal torus of $\GL_n$. We denote by
$x_i\in X(\BT_0)$ the $i$-th coordinate function and by
$\omega_i^\vee\in Y(\BT_0)$ the cocharacter that sends $a$ to the diagonal
matrix with coefficients $1$ at all positions except the $i$-th where the
coefficient is $a$. Let $\zeta$ be a generator of $\BF_{q^n}^\times$.
We have $N_w(\omega_i^\vee)=(\zeta,\zeta^q,\ldots,\zeta^{q^{n-1}})^{q^{n-i+1}}$.
Furthermore, $x_1$ induces an isomorphism $T\iso \BF_{q^n}^\times$ and
we identify these groups.
Let $s_i=(i,i+1)$. We have $\beta_i^\vee=\omega_1^\vee-\omega_{i+1}^\vee$ and
$N_w(\beta_i^\vee)=\zeta^{1-q^{n-i}}$. So,
$\mu_{m_i}=\langle \zeta^{1-q^{n-i}}\rangle$ by Proposition \ref{calculmi}.

Let $\theta\in\Irr(\BF_{q^n}^\times)$. We have
$c_\theta=\min\{a\ge 1|\theta(\zeta^{1-q^a})=1\}$ and
$n-1-|I_\theta|=|\{i\in\{1,\ldots,n-1\} | \theta(\zeta^{1-q^{n-i}})=1\}|$.
It follows that 
$n-1-|I_\theta|=\frac{n}{c_\theta}-1$.

Let $\psi$ be a linear regular character of $U_0$.
We deduce now from Corollary \ref{specE}
that $H^i(RE)^{\tilde{F}}=0$ for $i\not=0$, hence using
Corollary \ref{A centre connexe}, we obtain
$\Hom_{D^b(RT)}(RC,RC[i])^{\tilde{F}}=0$ for
$i\not=0$.

By Lemmas \ref{local} and \ref{EndK}, the canonical map
$b\CO G\to \End_{K^b(\CO T)}(C)^{\tilde{F}}$ is an isomorphism.
Lemma \ref{specE} shows that the corresponding assertion is true over $k$
as well.
\end{proof}

%\begin{proof}
%We know that $H^i(E)^{F'}=0$ for $i\not=0$ by Proposition \ref{descE}.
%Since $q^i\not=1$ in $k$ for $0<|i|\le r$ and since
%$H^i(kE)\simeq kH^i(E)$, it follows that $H^i(kE)=0$ for $i\not=0$.
%
%By Lemmas \ref{local} and \ref{EndK}, the canonical map
%$b\CO G\to \End_{K^b(\CO T)}(C)^{F'}$ is an isomorphism.
%Corollary \ref{specE} shows that the corresponding assertion is true over $k$
%as well.
%\end{proof}

\subsection{Brou\'e's conjecture}

\begin{thm}
\label{ADC}
Assume $\BG$ has type $A_{n-1}$ and $F$ is a split Frobenius endomorphism 
over $\BF_q$. Assume moreover that the order of $q$ modulo $\ell$ is equal to
$n$. 
Let $\CO Gb$ be the sum of blocks of $\CO G$ with non-zero defect. Let
$S$ be a Sylow $\ell$-subgroup of $T$. Then, $C_G(S)=T$.

The action of $\CO G$ on $\tilde{R}\Gamma_c(Y(\dot{w}),\CO)$ comes
from an action of $\CO Gb$ and the right action of $\CO T$ extends to
an action of $\CO N_G(T)$. The complex thus obtained induces a splendid
Rickard equivalence between $\CO Gb$ and $\CO N_G(T)$.
\end{thm}

\begin{proof}
The complex of $(\CO G,\CO T)$-bimodules
$C'$ (cf \S \ref{seclocal}) is $w$-stable. It follows that
it extends to a complex $D$ of
$(\CO G,\CO (T\rtimes\langle w\rangle))$-bimodules.

On the other hand, the complex $kC$ is isomorphic in
$D^b(k(G\times (T\rtimes F)^\oppose)\mMOD)$ to a bounded complex
$D'$ whose terms are finite dimensional.
There is a positive integer $N$ such that
$F^N$ acts trivially on $D'$. We take for $\nu$ a positive integer such that
$\ell^\nu\ge |k|$, $\ell^\nu\equiv 1 \pmod n$,
and $\nu\ge \nu_l(N)$ and we put $\tilde{F}=F^{\ell^\nu}$ as above.
There is
a positive integer $t$ prime to $\ell$ such that $\tilde{F}^t$ acts
trivially on $D'$. Let $e$ be a block idempotent of $b\CO G$. Then,
it follows from Proposition \ref{EndLambda} that $\End_{D^b(kG\otimes
(kT\rtimes \langle \tilde{F}\rangle)^\oppose)}(eD')\simeq Z(ekG)$, hence
$eD'$ is isomorphic
in $D^b(kG\otimes (kT\rtimes \langle \tilde{F}\rangle)^\oppose)$ to an 
indecomposable complex $X$.
The complexes $kX$ and $ekD$ have quasi-isomorphic restriction to
$kG\otimes kT$, hence (Clifford theory) they differ by tensor product
by a one-dimensional $k\langle \tilde{F}\rangle$-module $L$, where we
identify $w$ with $\tilde{F}$. In particular,
the canonical map
$$ekG\to \End_{k(T\rtimes w)^\oppose}^\bullet(ekD)$$
is an isomorphism. It follows from
\cite[Theorem 2.1]{Ri2} that $ekD$ is a two-sided tilting complex, \ie,
the canonical map
$\br(e)k(T\rtimes w)\to \End_{ekG}^\bullet(ekD)$ is an isomorphism as well.
We deduce that $D$ is a two-sided tilting complex for $b\CO G$ and
$\CO(T\rtimes w)$ (cf \eg \cite[Proof of Theorem 5.2]{Ri2}).

Let $Q$ be an $\ell$-subgroup of $G\times T^\oppose$.
We have $Y(\dot{w})^Q=\emptyset$ unless $Q$ is conjugate to a subgroup
of $\Delta S$ and $Y(\dot{w})^Q=T$ if $Q$ is a non-trivial subgroup
of $\Delta S$. It follows (cf \cite[proof of Theorem 5.6]{Rou1}) that 
$D$ is a splendid Rickard complex.
\end{proof}

\begin{rem}
For groups not untwisted of type $A$, the $\CO G$-module
$M=\Ind_{D(\BU_0)^F}^G\CO$ is not a progenerator in general.
It might nevertheless be possible to prove that
$R\Gamma_c(D(\BU_0)^F\backslash Y(\dot{w}),\CO)$ induces a derived
equivalence between $\End_{\CO G}(bM)$ and a quotient of
$\CO N_G(T)$.
\end{rem}

\end{document}